\magnification=1200
\input amstex
\documentstyle{amsppt}

\pageheight {49pc}
\vcorrection{-2pc}


\def\leftheadline {\line{\hfill\eightrm MARSHALL A. WHITTLESEY\hfill\folio}}
\def\rightheadline{\line{\hfill \eightrm FOLIATION BY GRAPHS OF CR MAPPINGS
\hfill \folio}}

\headline{\ifodd \pageno\rightheadline \else\leftheadline\fi}
\footline{}
\def\C{\hbox{\bf C$\,$}}

\def\Re{\hbox{\rm Re}}

\def\Nu{\hbox{H}}
\def\det{\hbox{\rm det\,}}
\def\M{\hbox{$\Cal M$}}

\def\i{\hbox{\bf i}}

\topmatter
\title  Foliation by graphs of CR mappings and a nonlinear Riemann-Hilbert
problem for smoothly bounded domains
\endtitle
\author Marshall A. Whittlesey
\endauthor
\affil University of California, San Diego
\endaffil
\address Department of Mathematics 0112, University of California, San Diego,
9500 Gilman Drive, La Jolla CA 92093-0112
\endaddress
\email mwhittlesey\@math.ucsd.edu
\endemail

\abstract 
Let $S$ be a generic $C^\infty $ CR manifold in $\Bbb C^\ell$,
$\ell\geq 2$, and let $M$ be a generic $C^\infty$ CR submanifold
of $S\times\Bbb C^m$, $m\geq 1$.  We prescribe conditions on $M$
so that it is the disjoint union of graphs of CR maps
$f:S\rightarrow\Bbb C^m$.  We also consider the special case
where $S$ is the boundary of a bounded, smoothly bounded
open set $D$ in $\Bbb C^\ell$, $\ell\geq 2$.  In this case we
obtain conditions which guarantee that such an $f$ above extends
to be analytic on $D$.  This provides a solution to a 
particular nonlinear Riemann-Hilbert boundary value problem
for analytic functions.

\endabstract
\date November 23, 2000 \enddate
\endtopmatter

\document
\NoBlackBoxes

\hsize=6.75truein
\hoffset=-.15in
\voffset=0.3 true in
\vsize=8.7truein
\nologo
Let $D$ be a bounded, smoothly bounded domain in $\Bbb C^\ell$,
$\ell\geq 2$, and let $M$ be a real $C^\infty$ submanifold of
$\partial D\times\Bbb C^m$, $m\geq 1$.  We here address the
question
of when there exists a mapping $f$ such that
$$\eqalign{(RH)&\left[\vbox{\hsize=4.5 true in \noindent
$f:\overline D\rightarrow\Bbb C^m$ is continuous on
$\overline D$ and analytic on $D$ such that the graph of $f$ over
$\partial D$ is contained in $M$.\vskip -0.1 in}\right.}$$

A problem where one is required to find an $f$ satisfying (RH) is
often called a Riemann-Hilbert problem; Riemann proposed such a
question for $\ell=m=1$ in 1851.  We shall refer to the problem of
finding an $f$ satisfying (RH) as the Riemann-Hilbert problem for
$M$. If an $f$ exists satisfying (RH) then we shall say that the
Riemann-Hilbert problem for $M$ is {\it solvable} and that $f$ is a
{\it solution}.  For $z\in\partial D$, let $M_z\equiv\{w\in\Bbb
C^m:(z,w)\in M\}$.  We will say here that the Riemann-Hilbert problem
(RH) is {\it linear} if for every $z\in\partial D$, $M_z$ is a real
affine subspace of $\Bbb C^m$.  We shall say that the Riemann-Hilbert
problem (RH) is {\it nonlinear} if it is not linear.

For the case $\ell=1$ we mention a few
references:[V,Sh1,Sh2,Fo,S,HMa,We1-5,B1,B2].  See [We4] for a useful
survey and reference list.  For $\ell\geq 2$, see [B1,B3,BD,D1,D2].

We will first address the following more general question.  Let $S$ be a
$C^\infty $ CR manifold in $\Bbb C^\ell$ (e.g., a real hypersurface in
$\Bbb C^\ell$) and let $M$ be a real $C^\infty$ submanifold of
$S\times
\Bbb C^m$, $m\geq 1$.  Let $z^0\in S$ and $U$ a neighborhood of $z^0$
in $S$.  Does there exist a CR map $f:U\rightarrow \Bbb C^m$ whose
graph in $U\times
\Bbb C^m$ is contained in $M$?  We shall establish conditions under which
the answer to this question is yes.  For the case where $S$ is a
complex manifold, this question is addressed by theorems in
[Fr2-4,Kr,So1-2].  We consider the case of more general $S$.  Applying
our result to the case where $S$ is the boundary of an open set in
$\Bbb C^\ell$ with $C^\infty$ boundary (so $S$ is a real hypersurface
in $\Bbb C^\ell$,) we shall establish conditions where the
Riemann-Hilbert problem for $M$ is solvable.

For the more general question of the previous paragraph, assumptions
we shall make will imply that the set $M$ is a CR manifold.  Under
conditions outlined in Theorem 1, we will prove the existence of a
special involutive subbundle of the real tangent bundle to $M$ which
possesses certain null properties relative to the Levi form of
$M$. The Frobenius theorem will guarantee the existence of integral
submanifolds of this bundle.  Locally, these integral submanifolds
will turn out to be graphs of CR maps from open subsets of $S$ to
$\Bbb C^m$.  In Theorem 2 and Theorem 3 we establish conditions which
guarantee the existence of graphs of CR maps contained in $M$ which
are defined on all of $S$.  In Theorem 4 we show that under some
conditions these graphs possess a particular extremal property.  In
Corollary 1 and Corollary 2, we will assume that $S$ is the boundary
of a bounded, $C^\infty$-bounded open set in $\Bbb C^\ell$ and
establish conditions where the CR maps found in Theorem 2 and Theorem
3 extend to be analytic on $D$.  These maps are solutions to the
Riemann-Hilbert problem for $M$.  The conditions we establish for the
solvability of (RH) are not always necessary: when they are satisfied,
we obtain graphs with strong properties which we shall describe.

We are grateful to Professors Salah Baouendi, Linda Rothschild and
Peter Ebenfelt for useful conversations on this topic.

In studying subsets of $\Bbb C^\ell\times\Bbb C^m$ we shall generally
label points in $\Bbb C^\ell$ with $z=(z_1,z_2,...,z_\ell)$ and those
in $\Bbb C^m$ with $w=(w_1,w_2,...,w_m)$.  If $P(z,w)$ is a function
on an open subset of $\Bbb C^\ell\times\Bbb C^m$ then let $\partial P$
denote the $1$-form $\sum_{i=1}^\ell{\partial P\over\partial
z_i}dz_i+\sum_{i=1}^m {\partial P\over\partial w_i}dw_i$ and let
$\partial_wP$ denote the $1$-form $\sum_{i=1}^m{\partial
P\over\partial w_i}dw_i$.  Let $[X,Y]$ denote the Lie bracket of two
vector fields $X,Y$, $\langle A,B\rangle$ denote the canonical action
of $A$ on $B$ where for some $n$, $A$ is an $n$-cotangent (i.e., an
element of the $n^{th}$ exterior algebra of the cotangent space) at a
point on a manifold, $B$ an $n$-tangent at that point; see [Bo,
p. 10].  If for $i=1$ to $n$, $\phi_i$ is a cotangent at a point and
$\psi_i$ is a tangent at the same point then
$\langle\phi_1\wedge\phi_2\wedge...\wedge
\phi_n,\psi_1\wedge\psi_2\wedge...\wedge
\psi_n\rangle$ is the determinant of the matrix whose $i,j$ component
is $\langle \phi_i,\psi_j\rangle$.  We use the same notation to denote
the action of $A$ on $B$ where $A$ is an $n$-form on a manifold and
$B$ a section of the $n^{th}$ exterior algebra of the tangent bundle,
the action being defined pointwise.  If $B$ is a vector bundle on an
open subset $A$ of a $C^\infty$ manifold then $\Gamma(A,B)$ will
denote the set of $C^\infty $ sections of $B$ over $A$.  We shall also
make use of the notion of {\it generic} CR manifold (see [BER, p. 9])
and CR manifold of (Bloom-Graham-Kohn) {\it finite type} $\tau$ (see
[BER, pp. 17-18]).

If $\M$ is a $C^\infty$ generic CR submanifold of $\Bbb C^n$, then let
$T\M$ denote the real tangent bundle to $\M$, $\Bbb CT\M$ the
complexified tangent bundle to $\M$, $\Bbb CT_x\M$ the fiber of $\Bbb
CT\M$ over $x\in \M$, $\Bbb CT^{1,0} \M$ be the $(1,0)$ subbundle of
$\Bbb CT\M$ and $\Bbb CT_x^{1,0} \M$ the fiber of $\Bbb CT^{1,0} \M$
over $x$.  (We use similar notation for the $(0,1)$ subbundle.)  We
shall make use of the Levi map of $\M$, and we follow the definition
in [BER]: Let $P_x:\Bbb CT_x\M\rightarrow
\Bbb CT_x\M/(\Bbb CT_x^{1,0} \M\oplus\Bbb CT_x^{0,1} \M)$ be projection
and let the Levi map at $x\in \M$ be
$$\eqalign{{\Cal L}_x:\Bbb CT_x^{1,0} \M\times \Bbb CT_x^{1,0}
\M&\rightarrow\Bbb CT_x\M/(\Bbb CT_x^{1,0} \M\oplus\Bbb CT_x^{0,1} \M)
\cr {\Cal
L}_x(X_x,Y_x)&:={1\over 2i}P_x([X,\overline Y](x)),\cr}$$ where $X,Y$
are $C^\infty$ $(1,0)$ vector fields on $\M$ near $x$ such that
$X(x)=X_x$ and $Y(x)=Y_x$.  As is well known, the Levi map at $x$ does
not depend on these smooth extensions $X,Y$; it depends only on $X_x$
and $Y_x$, so we get a smooth bundle map ${\Cal L}:\Bbb CT^{1,0}
\M\times \Bbb CT^{1,0}
\M  \rightarrow\Bbb
CT\M/( \Bbb CT^{1,0} \M\oplus \Bbb CT^{0,1}
\M)$.  Suppose that $\M$ has defining functions $\phi_i$, $i=1,2,...,d$.
Note that ${\Cal L}_x(X_x,Y_x)=0$
if and only if $\langle\partial \phi_i, [X,\overline Y]\rangle(x)=0$
for $i=1$ to $d$ and for $X,Y$ as before.  By Cartan's identity (see
[Bo,p. 14, Lemma 3]) these equations are equivalent to the equations
$\langle\overline\partial\partial \phi_i,X\wedge
\overline Y\rangle(x)=0$ for $i=1$ to $d$.

We follow the definition of CR function given in [Bo]: a CR
function on a $\M$ is a $C^1$ function $u$ such that for every section
$X\in\Gamma(\M,\Bbb CT^{0,1}\M)$ we have $Xu=0$. 

We assume that $S$ satisfies the following properties.

$$\eqalign{(1)&\left[\vbox{\hsize=4.5 true in Let $S$ be a generic
$C^\infty $ CR manifold in $\Bbb C^\ell$, $\ell\geq 2$, of CR
codimension $c<\ell$.  Specifically, we assume that $S$ has $C^\infty$
real valued defining functions $p_1,p_2,...,p_c$ defined in a
neighborhood $\Cal{N}$ of $S$ (i.e., for $i=1$ to $c$, $p_i:{\Cal
N}\rightarrow \Bbb R$ is $C^\infty$, $S=\{z\in{\Cal
N}:p_i(z)=0\,\,\hbox{for all $i=1$ to $c$}\}$ and the $1$-forms
$\{\partial p_i(z)\}_{i=1}^c$ are linearly independent $1$-forms over
$\Bbb C$ for all $z\in S$.)  Let $\Nu^S=\Bbb CT^{1,0}S$ be the $(1,0)$
tangent bundle to $S$ with fibers $\Nu^S_z$, i.e. $$\Nu^S_z=\{
\sum_{i=1}^\ell a_i{\partial\over\partial z_i}\in\Bbb CT_zS:a_i\in\Bbb C,
i=1,2,...,c\}.$$ 
\vskip -1.25 in
}
\right .}$$
Later we shall assume that $S$ is the $C^\infty $ boundary of 
a bounded domain $D$.

We assume $M$ satisfies the following properties.

$$\eqalign{(2)&\left[\vbox{\hsize=4.5 true in Let $d$ be an integer,
$1\leq d\leq m$, and let $M$ be a $C^\infty $ CR submanifold of
$S\times\Bbb C^m$ of CR codimension $c+d$ in $\Bbb C^{\ell+m}$ such
that if we write $M_z\equiv\{w\in\Bbb C^m:(z,w)\in M\}$ then for all
$z\in S$, $M_z$ is a nonempty, generic, Levi nondegenerate CR
submanifold of CR codimension $d$ in $\Bbb C^m$.  Specifically, let
${\Cal U}\subset\Bbb C^{\ell+m}$ be an open set meeting $S\times
\Bbb C^m$, $q_i: {\Cal U} \rightarrow \Bbb R$ 
be $C^\infty$ for $i=1,2,...,d$ such that for any $(z^0,w^0)\in {\Cal
U}$ the $1$-forms $\{\partial_wq_i(z^0, w^0)\}$, $i=1$ to $d$, are
linearly independent over $\Bbb C$.  Let $M = \{(z, w)
\in {\Cal U}\cap (S\times\Bbb C^m)|q_1(z, w) = 
q_2(z,w)=...=q_d(z,w)=0\}$ and $\Nu^M\equiv\Bbb CT^{1,0}M$ be the
$(1,0)$ tangent bundle to $M$ with fibers $\Nu^M_{(z,w)}$ for
$(z,w)\in M$.
\vskip -1.07 in
}
\right .}$$
Note that since $M_z$ is generic, we must have $d\leq m$.
The following correspondence is convenient:
$$\Nu^M_{(z, w)}  \cong$$ $$\left\{\sum\limits^\ell_{i=1}a_i{\partial \over
\partial z_i}  +
\sum\limits^m_{i=1} b_i{\partial \over \partial w_i}\in \Nu^S_z 
\oplus \Bbb CT_w^
{1,0}\Bbb C^m\left|\sum\limits^\ell_{i=1}\right. a_i{\partial q_k\over
\partial z_i} (z,w)  +
\sum\limits^m_{i=1} b_i{\partial q_k\over \partial w_i} (z,w)
= 0, k=1..d\right\}.$$ 
Note that $\{p_i\}_{i=1}^c\cup\{q_i\}_{i=1}^d$
is a set of defining functions for $M$.  Let $\pi$ be projection from
$\Bbb CT\Bbb C^{\ell+m}$ to $\Bbb CT
\Bbb C^\ell$, or the restriction to the corresponding mapping from
$\Nu^M$ to $\Nu^S$.  We will use $\pi_{(z,w)}$ to denote projection on the
individual fibers: $\pi_{(z,w)}:\Bbb CT_{(z,w)}\Bbb C^{\ell+m}\rightarrow
\Bbb CT_z\Bbb C^\ell$ and $\pi_{(z,w)}:\Nu^M_{(z,w)}\rightarrow\Nu^S_z$.

We shall investigate when there exists a $C^\infty$ CR map $f:S\rightarrow
\Bbb C^m$ such that $q_i(z,f(z))=0$ for $z\in S$ and $i=1$ to $d$.  

Let $$V(z, w)=\{v_{zw} \in
\Nu^M_{(z, w)}:\pi_{(z,w)}(v_{zw})=0\};\eqno(3)$$ 
that is, if $v_{zw}\in V(z,w)$ then $v_{zw}$ has no terms involving
${\partial\over\partial z_i}$, $i=1,2,...,\ell$.  (This is the space
of ``vertical" $(1,0)$ tangents to $M$ at $(z, w)$.)  By (2), $V(z,w)$
has dimension $m-d$ for all $(z,w)\in M$.  For any $(z^0,w^0)\in M$,
$V(z^0,w^0)=\{v_{z^0w^0}\in\Bbb CT_{(z^0,w^0)}\Bbb
C^{\ell+m}:\pi_{(z^0,w^0)}(v_{z^0w^0})=0\,\hbox{and}\,
\langle\partial_wq_i(z^0,w^0),v_{z^0w^0}\rangle=0,i=1,2,...,d\}$.  Thus
$V(z,w)$ is calculated by solving a system of linear equations whose
coefficients are $C^\infty$ functions of $(z,w)$ and whose rank is
constant for $(z,w)\in M$ near $(z^0,w^0)$.  We may thus calculate
$m-d$ $C^\infty$ vector fields $v^i,i=1,2,...,m-d$, near any point
$(z^0,w^0)\in M$ such that near $(z^0,w^0)$,
$\{v^i(z,w)\}_{i=1}^{m-d}$ forms a basis for $V(z,w)$.  Thus we obtain
a complex $(m-d)$-dimensional $C^\infty $ $(1,0)$ vector bundle $V$
over $M$.  Let $\overline V$ denote the bundle whose fiber at $(z,w)$
is $\overline {V(z,w)}$.  Let ${\Cal L}^M$ denote the Levi map for
$M$.  Next let $$N(z, w)=\{n_{zw}\in
\Nu^M_{(z,w)}:{\Cal L}^M_{(z,w)}(n_{zw}, v_{zw})=0\,\forall v_{zw}\in 
V(z,w)\}.\eqno(4)$$
Because $V$ is a bundle, any $v_{zw}$ is the value at $(z,w)$ of some
element of $\Gamma(M,V)$, so in (4) it is equivalent to demand that
${\Cal L}^M_{(z,w)}(n_{zw}, v(z,w))=0$ for all $v\in\Gamma(M,V)$.

Suppose we have fixed $z^0\in S$ and we choose $C^\infty $ vector
fields $s_i$, $i=1$ to $\ell-c$ defined on an open set $G$ containing
$z^0$, such that for $z\in G$,
$$\{s_i(z)\}_{i=1}^{\ell-c}\eqno(5)$$ is a basis for $\Nu^S_z$. 
Let $$G^M\equiv M\cap(G\times\Bbb C^m).\eqno(6)$$

\proclaim{Proposition 1}
Assume that $S$ and $M$ satisfy (1) and (2), respectively.  For any
$(z,w)\in M$ and any value of $c$, the projection $\pi_{(z,w)}$ of
$N(z,w)$ to $\Nu^S_z$ is injective.  Thus for every $(z, w)
\in M$, the complex dimension of $N(z, w)$ is less than or equal to
$\ell-c$.  If the complex dimension of $N(z,w)$ equals $\ell-c$ for
every $(z,w)\in G^M$ then $N$ is a complex $C^\infty$ vector bundle of
dimension $\ell-c$ over $G^M$ and $\pi_{(z,w)}:N(z,w)\rightarrow
\Nu^S_z$ is an isomorphism.  The complex dimension of $N(z,w)$ is
exactly $\ell-c$ in the following three cases: (i) when $d=1$, (ii)
when $d=m$, and (iii) when the zero sets of $q_i$, $i=2$ to $d$, are
Levi flat in $\Bbb C^{\ell+m}$ near $M$ (i.e. the Levi form of those
surfaces is totally degenerate.)
\endproclaim

\noindent{\bf Proof:}
If $n_{zw}\in N(z,w)$ satisfies the property that
$\pi_{(z,w)}(n_{zw})=0$ then $n_{zw}\in V(z,w)$.  Then $$\langle
\overline\partial\partial q_i(z,w),n_{zw}\wedge \overline v_{zw}
\rangle =0\eqno(7)$$ for all
$v_{zw}\in V(z,w)$ and $i=1$ to $d$.  Note that $V(z,w)$ may be
regarded as the $(1,0)$ tangent space at $w$ to $M_z$, which is Levi
nondegenerate by (2).  Then (7) implies that $n_{zw}$ is in the null
space of the Levi form of $M_z$, so $n_{zw}=0$.  This proves that
$\pi_{(z,w)}:N(z,w)\rightarrow \Nu^S_z$ is injective and the complex
dimension of $N(z,w)$ is less than or equal to the dimension of
$\Nu^S_z$, which is $\ell-c$.  Now fix an arbitrary $(z^0,w^0)\in M$
and fix $a^0=(a_1^0,a_2^0,...,a_\ell^0)\in\Bbb C^\ell$ such that
$\sum_{i=1}^\ell a_i^0({\partial p_j/\partial z_i})(z^0)$ $=0$ for
$j=1$ to $c$.  For $i=1$ to $\ell$ let $a_i$ be a $C^\infty$ complex
function defined on $S$ near $z^0$ such that $a_i(z^0)=a_i^0$ and
$\sum_{i=1}^\ell$ $a_i(z)({\partial p_j/\partial z_i})(z)=0$ for $j=1$
to $c$ and $z$ near $z^0\in S$ (i.e. $A(z)\equiv\sum_{i=1}^\ell a_i(z)
{\partial\over\partial z_i}\in
\Nu^S_z$.)  For $(z,w)\in M$ near $(z^0,w^0)$, necessary and sufficient
conditions for the existence of a
$\Bbb C^m$-valued mapping $b=(b_1,b_2,...,b_m)$ on $M$ such that
$B(z,w)\equiv\sum_{i=1}^\ell a_i(z)({\partial /\partial z_i})+
\sum_{i=1}^m b_i(z,w)({\partial /\partial w_i})$ belongs to $N(z,w)$
(and so $\pi_{(z,w)}(B(z,w))=A(z)$)
are

$$\sum\limits^\ell_{i=1} a_i(z) {\partial q_k\over \partial z_i}
(z, w)+\sum\limits^m_{i=1} b_i(z,w)
{\partial q_k\over \partial w_i}(z, w) = 0 \eqno(8)$$
for $k=1$ to $d$,
$$\langle \overline \partial\partial p_j, \left(
\sum\limits^\ell_{i=1} a_i
{\partial \over \partial z_i} + \sum\limits^m_{i=1} b_i
{\partial \over \partial w_i}\right)\wedge
\left(\sum\limits^m_{i=1} \overline v_i {\partial \over
\partial \overline w_i}\right)\rangle(z,w) =0 \eqno(9)$$
and
$$\langle \overline \partial\partial q_k, \left(
\sum\limits^\ell_{i=1} a_i
{\partial \over \partial z_i} + \sum\limits^m_{i=1} b_i
{\partial \over \partial w_i}\right)\wedge
\left(\sum\limits^m_{i=1} \overline v_i {\partial \over
\partial \overline w_i}\right)\rangle(z,w) =0 \eqno(10)$$
for $j=1$ to $c$, $k=1$ to $d$ and for all vertical $(1,0)$ vector
fields $\sum\limits^m_{i=1} v_i {\partial \over \partial
w_i}\in\Gamma(M,V)$.  Condition (9) is vacuous: $\overline
\partial\partial p_j$ contains only terms involving $d\overline
z_{i_1}\wedge dz_{i_2}$ but $\sum\limits^m_{i=1} \overline v_i
{\partial \over
\partial \overline w_i}$ contains no terms with ${\partial \over
\partial \overline z_i}$ or ${\partial \over
\partial z_i}$, so the left side of (9) is automatically zero.

Thus conditions (8) and (10) impose the requirement that
$b_1(z,w),b_2(z,w),..., b_m(z,w)$ must satisfy a nonhomogeneous system
of linear equations.  By (2), the dimension of $V(z,w)$ is $m-d$ for
all $(z,w)\in M$.  For a small open set $U\subset G^M$ containing
$(z^0,w^0)$, we may choose elements $v^j=\sum_{i=1}^m
v^j_i{\partial\over\partial w_i}\in\Gamma(U,V)$, $j=1$ to $m-d$, such
that $\{v^j(z,w)\}_{j=1}^{m-d}$ is a basis for $V(z,w)$ for all
$(z,w)\in U$.  For equations (8),(10) to hold for $(z,w)\in U$, it is
equivalent for the following system of equations to hold for $(z,w)\in
U$:
$$\sum\limits^\ell_{i=1} a_i(z) {\partial q_k\over \partial z_i}
(z, w)+\sum\limits^m_{i=1} b_i(z,w)
{\partial q_k\over \partial w_i}(z, w) = 0 \eqno(11)$$
for $k=1$ to $d$
and
$$\langle \overline \partial\partial q_k, \left(
\sum\limits^\ell_{i=1} a_i
{\partial \over \partial z_i} + \sum\limits^m_{i=1} b_i
{\partial \over \partial w_i}\right)\wedge
\left(\sum\limits^m_{i=1} \overline v^j_i {\partial \over
\partial \overline w_i}\right)\rangle(z,w) =0 \eqno(12)$$
for $k=1$ to $d$ and $j=1$ to $m-d$.

Our system has $d$ equations from (11) and $d(m-d)$ equations from (12)
(not necessarily independent), so is one of $d+d(m-d)$ equations in
$m$ unknowns $b_1(z,w),...,b_m(z,w)$.  (Note that if $d=1$ or $d=m$
then $d+d(m-d)=m$.)  

Suppose that $N(z,w)$ has complex dimension $\ell-c$ for all $(z,w)\in
G^M$.  Since $\pi_{(z,w)}:N(z,w)\rightarrow \Nu^S_z$ is an injection
and $\Nu^S_z$ also has complex dimension $\ell-c$,
$\pi_{(z,w)}:N(z,w)\rightarrow \Nu^S_z$ is an isomorphism.  Thus given
any element $A\in\Gamma(G,\Nu^S)$ there exists a unique $B\in
\Gamma(G^M,\Nu^M)$ such that $\pi_{(z,w)}(B(z,w))=A(z)$.  That means
that for $(z,w)\in G^M$, the system (11),(12) has precisely one
solution for $b(z,w)$ given a fixed $A(z)$, so there exists a
subsystem of $m$ linear equations in the $b_j(z,w)$ whose solution is
the same as for (11),(12), for $(z,w)$ near some point $(z^1,w^1)\in
G^M$.  That implies that the $b_j$ are all $C^\infty$ functions near
$(z^1,w^1)$, since they are the unique solution to a system of $m$
equations in $m$ unknowns with coefficients which are $C^\infty$ in
$(z,w)$.  This holds near an arbitrary $(z^1,w^1)\in G^M$, so $B$ is
$C^\infty$ (and we write $B\in\Gamma(G^M,\Nu^M)$.)  Let $n_i$ be the
unique $C^\infty$ $(1,0)$ vector field defined on $G^M$ such that
$n_i(z,w)\in N(z,w)$ for $(z,w)\in G^M$ and
$$\pi_{(z,w)}(n_i(z,w))=s_i(z).\eqno(13)$$ Then $\{n_i(z,w)\}$ is a
basis for $N(z,w)$ for all $(z,w)\in G^M$ since
$\pi_{(z,w)}:N(z,w)\rightarrow \Nu^S_z$ is an isomorphism.  We
conclude that the $N(z,w)$ form a $C^\infty $ bundle over $G^M$; since
sets such as $G^M$ cover $M$, we have a bundle on all of $M$ which we
call $N$.

The system of equations associated to (11),(12) which is
homogeneous in the $b_i(z,w)$ (i.e., $A(z)=0$) is
$$\sum\limits^m_{i=1} b_i(z,w) {\partial q_k\over \partial w_i}(z, w)
= 0 \eqno(11\hbox{$'$})$$ for $k=1$ to $d$ and
$$\langle \overline \partial\partial q_k, \left(
\sum\limits^m_{i=1} b_i
{\partial \over \partial w_i}\right)\wedge
\left(\sum\limits^m_{i=1} \overline v^j_i {\partial \over
\partial \overline w_i}\right)\rangle(z,w) =0 \eqno(12\hbox{$'$})$$
for $k=1$ to $d$ and $j=1$ to $m-d$.
As noted earlier, $b(z,w)=0$ is the unique solution to the homogeneous
system because $M_z$ is Levi nondegenerate for all $z\in S$.  If $d=1$
or $d=m$ then the number of equations in the $b_j$ is $d+d(m-d)=m$, so
then the nonhomogeneous system (11),(12) has exactly one solution for
the $b_j(z,w)$, $j=1,2,...,m$.  Then what we have just shown is that
for $d=1$ or $m$ and $(z,w)\in M$, every element
$A(z)$ in $\Nu^S_z$ is the image under $\pi_{(z,w)}$ of exactly
one element $B(z,w)$ in $N(z,w)$; i.e. the projection $\pi_{(z,w)}:
N(z,w)\rightarrow
\Nu^S_z$ is bijective and the complex dimension of $N(z,w)$ is equal 
to $\ell-c$.

To prove (iii), we must determine the number of independent equations
in the $b_i$ arising from (11) and (12).  From (11) we obtain $d$ of
them (one for each $q_i$).  From (12) we obtain $m-d$ equations
arising from the second order equations involving $q_1$, since the
rank of $V$ is $m-d$; all other second order equations are vacuous by
the Levi flatness associated with the other $q_i$.  Thus we have a
total of $m$ equations in $(b_1,b_2,...,b_m)$; since (as before) the
associated homogeneous system (11\hbox{$'$}),(12\hbox{$'$}) has exactly
one solution (by the Levi nondegeneracy of $M_z$ again), the
nonhomogeneous system (11),(12) does also.  Following the end of the
argument above, we conclude that the dimension of $N(z,w)$ is $\ell-c$
for all $(z,w)\in M$.$\hfill\square$

We use $\overline N$ to denote the bundle with fibers $\overline
{N(z,w)}$.  If $B$ is a vector bundle over a manifold $A$, then we say
that a set of vector fields $\{L_i\}_{i\in I}$ is a {\it local basis}
for $B$ near a point in $A$ if for all $x$ in $A$ near that point
$\{L_i(x)\}_{i\in I}$ is a basis for the fiber $B_x$ of $B$ over $x$.

Proposition 2 will establish conditions where the spaces $N(z,w)$
together compose an involutive bundle over $M$.  However, we first
need the following lemma.  We say that a vector field is a commutator
of length $\sigma\geq 2$ if it has the form
$[Y_1,[Y_2,[Y_3,...,[Y_{\sigma-1},Y_\sigma]]...]$, for vector fields
$Y_i$, $i=1,2,3,...\sigma$.

\proclaim{Lemma 1}  Let $T$ be any commutator of the vector fields $n_1,n_2,
...,n_{\ell-c},\overline n_1,\overline n_2,...,\overline n_{\ell-c}$.
Then for all $v_1,v_2\in\Gamma(G^M,V)$, $$[T,v_1+\overline v_2]
\in\Gamma(G^M,V\oplus\overline V).\eqno(14)$$
(Furthermore, (14) holds if $T$ is a linear combination of commutators
of the $n_i,\overline n_i$.)
\endproclaim

\noindent{\bf Proof:} The parenthetical sentence of Lemma
1 follows from the second sentence because (14) is linear in $T$.  We
prove the second sentence by induction on the length of $T$.  If the
length of $T$ is $1$ then $T=n_i$ or $T=\overline n_i$ for some $i=1$
to $\ell-c$.  If $T=n_i$,
$$[n_i,v_1+\overline v_2] =[n_i,v_1]+[n_i,\overline v_2].\eqno(15)$$
The first term on the right hand side of (15) belongs to
$\Gamma(G^M,\Nu^M)$ by involutivity of $\Nu^M$.  Since from (13) the
coefficients of $n_i$ in ${\partial\over\partial z_j}$ depend only on
$z\in S$ for $j=1$ to $\ell$, $v_1$ annihilates these coefficients.
Thus $[n_i,v_1]$ has no terms involving ${\partial \over\partial
z_j}$, so $\pi_{(z,w)}([n_i,v_1](z,w))=0$ and
$[n_i,v_1]\in\Gamma(G^M,V)$ by definition of $V$.  The second term on
the right hand side of (15) belongs to
$\Gamma(G^M,\Nu^M\oplus\overline \Nu^M)$ by definition of $N$.  Write
$[n_i,\overline v_2]=h^1+\overline h^2$ for $h^j\in\Gamma(G^M,\Nu^M)$,
$j=1,2$.  For the same reason as with the first term of (15),
$\pi_{(z,w)}([n_i,\overline v_2](z,w))=0$, i.e.
$\pi_{(z,w)}(h^1(z,w)+\overline h^2(z,w))=0$.  Thus
$\pi_{(z,w)}(h^1(z,w))$ and $\pi_{(z,w)}(\overline h^2(z,w))$ are
negatives of each other, but the former belongs to $\Nu^S_z$ and the
latter to $\overline\Nu^S_z$.  Since these two spaces meet only in
$\{0\}$, $\pi_{(z,w)}(h^i(z,w))=0$ for $i=1,2$.  Thus $h^1(z,w)$ and
$h^2(z,w)$ belong to $V(z,w)$, so $[n_i,\overline
v_2]\in\Gamma(G^M,V\oplus\overline V)$, as desired.  We conclude that
the right hand side of (15) belongs to $\Gamma(G^M,V\oplus\overline
V)$.  If $T=\overline n_i$ then $[T,v_1+\overline v_2]=[\overline n_i,
v_1+\overline v_2]=\overline{[n_i,\overline v_1+v_2]}$ which belongs
to $\Gamma(G^M,V\oplus\overline V)$ by what we just showed in the case
when $T=n_i$.  This proves the theorem if the length of $T$ is $1$.

For the remainder of the proof, the commutators referred to as $T,T',
T''$ will be commutators of vector fields in the set $\{n_i,\overline
n_i:i=1,2,...,\ell-c\}$.  Now suppose that Lemma 1 is true for all
commutators $T$ of length less than or equal to $\lambda$.  Then we
must show that if $T$ is a commutator of length $\lambda+1$, it
satisfies (14).  Then $T=[n_i,T']$ or $T=[\overline n_i,T']$ for some
commutator $T'$ of length $\lambda$.  We claim that it suffices to
prove (14) for all $T$ of the form $[n_i,T']$ where $T'$ is a
commutator of length $\lambda$; once this is done, if we write
$T''=[\overline n_i,T']$, then $[T'',v_1+\overline v_2]= [[\overline
n_i,T'],v_1+\overline v_2]=\overline{[[n_i,\overline T'],
\overline v_1+v_2]}$.  The expression $[[n_i,\overline T'],
\overline v_1+v_2]$ belongs to $\Gamma(G^M,V\oplus\overline V)$
since $\overline T'$ is a commutator of length $\lambda$.  Thus
$[T'',v_1+\overline v_2]=\overline{[[n_i,\overline T'],
\overline v_1+v_2]}$ also belongs to $\Gamma(G^M,V\oplus\overline V)$.

Now assume that $T=[n_i,T']$ for some commutator $T'$ of length
$\lambda$.  Then, using Jacobi's identity,
$$\eqalign{[T,v_1+\overline v_2]=[[n_i,T'],v_1+\overline v_2]
&=-[[T',v_1+\overline v_2],n_i]-[[v_1+\overline v_2,n_i],T']\cr
&=[n_i,[T',v_1+\overline v_2]]-[T',[n_i,v_1+\overline
v_2]].\cr}\eqno(16)$$ The vector fields $[T',v_1+\overline
v_2]$ and $[n_i,v_1+\overline v_2]$ both belong to
$\Gamma(G^M,V\oplus\overline V)$ by the induction hypothesis.  Then
(16) does as well, again by applying the induction hypothesis to each
term of (16).  That implies that $[T,v_1+\overline
v_2]\in\Gamma(G^M,V\oplus\overline V)$ also, so $T$ satisfies (14).
By induction, the second sentence of Lemma 1 is proven.
$\hfill\square$

\proclaim{Proposition 2}
Assume that $S$ and $M$ satisfy (1) and (2), respectively.  If
$N(z,w)$ has dimension $\ell-c$ for all $(z,w)\in M$, then $N$ is an
involutive $C^\infty$ subbundle of $\Nu^M$.
\endproclaim

\noindent{\bf Proof:}
$N$ is a $C^\infty$ bundle from Proposition 1.  To see that $N$ is
involutive, it suffices to verify this in a neighborhood of each point
$(z^0,w^0)\in M$.  Furthermore, it suffices to check the condition of
involutivity on the local basis $\{n_i\}$ from (13): we must show that
$[n_i,n_j](z,w)\in N(z,w)$ for all $(z,w)$ near $(z^0,w^0)$ in $M$.
Let $v_{(z^0,w^0)}\in V(z^0,w^0)$, $v\in\Gamma(M,V)$ such that
$v(z^0,w^0)=v_{(z^0,w^0)}$.  Then by Lemma 1, $[[n_i, n_j],\overline
v]\in \Gamma(G^M,V\oplus\overline V)
\subset\Gamma(G^M,\Nu^M\oplus\overline\Nu^M)$.
We already know that $\Nu^M$ is involutive, so $[n_i,n_j](z,w)\in
\Nu^M_{(z,w)}$.
Thus ${\Cal L}^M_{(z,w)}( [n_i,n_j](z,w),v_{(z,w)})=0$ for all
$v_{(z,w)}\in V(z,w)$, so by the
definition of $N$, $[n_i,n_j](z,w)\in N(z,w)$.  Thus $N$ is
involutive.\hfill $\square$
\bigskip

Note that Proposition 1 implies that $V(z,w)\cap N(z,w)=\{0\}$.  If
the dimension of $N(z,w)$ is $\ell-c$, it also implies that
$V(z,w)\oplus N(z,w)=\Nu^M_{(z,w)}$: we know that the dimension of
$V(z,w)$ is $m-d$ and the dimension of $\Nu^M_{(z,w)}$ is
$(\ell+m)-(c+d)$.  Thus the dimension of $\Nu^M_{(z,w)}$ is the sum of
the dimensions of $V(z,w)$ and $N(z,w)$; since $V(z,w)\cap
N(z,w)=\{0\}$, we must have $V(z,w)+N(z,w)=\Nu^M_{(z,w)}$ also.

Our intention is to construct a map $f:S\rightarrow\Bbb C^m$ whose
graph lies in $M$ and passes through a point $(z^0,w^0)\in M$, such
that the $(1,0)$ tangent space to the graph of $f$ at $(z,f(z))$ is
$N(z,f(z))$.  This will imply that $f$ is a CR mapping on $S$; see the
proof of Theorem 1.  We say that a real $C^\infty$ manifold $H$ of
real dimension $a$ is {\it foliated} by a class of submanifolds ${\Cal
C}$ of real dimension $b$ near a point $y\in H$ 
if there exists a $C^\infty$ diffeomorphism from
$Q\equiv\{(x_1,x_2,...,x_a)\in
\Bbb R^a: 0<x_i<1,i=1,2,...,a\}$ onto a neighborhood of $y$ in $H$
such that for any constants $c_i$ between $0$ and $1$,
$i=b+1,b+2,...,a$, the image of the set $\{(x_1,x_2,...,x_a)\in
Q:x_i=c_i,i=b+1,b+2,...,a\}$ is an open subset of some member of ${\Cal
C}$.  We shall say that $H$ is foliated by the manifolds in ${\Cal C}$
if $H$ is foliated by ${\Cal C}$ near each of its points.

Lemma 2 will be needed in part of the proof of Theorem 1.

\proclaim{Lemma 2} Suppose that $\phi_1,\phi_2,...,\phi_c,\phi_{c+1}$
are elements of the cotangent space $\Bbb CT^*_0(\Bbb R^c\times\Bbb
C^k)$ and that $\langle\phi_i,T\rangle=0$ for all tangent vectors
$T\in\Bbb CT_0(\{0\}\times\Bbb C^k)$.  Then
$\phi_1\wedge\phi_2\wedge\phi_3\wedge...\wedge\phi_{c+1}=0$.
\endproclaim

\noindent{\bf Proof:} Assume that $\Bbb R^c$ has
coordinates $x_1,x_2,...,x_c$ and $\Bbb C^k$ has coordinates
$\zeta_1,\zeta_2,..., \zeta_k$.  Then we may write
$\phi_j=\sum_{i=1}^c\alpha^j_idx_i+\sum_{i=1}^k\beta^j_id\zeta_i+
\sum_{i=1}^k\gamma^j_id\overline
\zeta_i$.  Since $\langle\phi_i,T\rangle=0$ for all tangent vectors
$T\in\Bbb CT_0(\{0\}\times\Bbb C^k)$, we have $\beta^j_i=\gamma^j_i=0$
for all $i=1$ to $k$ and $j=1$ to $c+1$.  Thus the
$\phi_j=\sum_{i=1}^c\alpha^j_idx_i$ may be regarded as elements of
$\Bbb CT^*_0\Bbb R^c$ for $j=1$ to $c+1$; taking the wedge product of
all $c+1$ of them is identically zero, since $c+1>c$.  (The wedge
product of linearly dependent vectors is zero, and any $c+1$ vectors
in a $c$-dimensional space are linearly dependent.) $\hfill\square$

\proclaim{Theorem 1}  Suppose that $S$, $M$ satisfy (1),(2), respectively,  
and that at every point of $S$, $S$ is of (Bloom-Graham-Kohn) finite
type $\tau$.  Then the following conditions are equivalent:

(I)  The dimension of $N(z,w)$ is $\ell-c$ for
all $(z,w)\in M$ and for every $k=1$ to $d$ we have that
$$\langle\partial p_1\wedge\partial p_2\wedge\partial p_3\wedge...\wedge
\partial p_c\wedge\partial q_k,F_1\wedge
F_2\wedge F_3\wedge...\wedge F_{c+1}\rangle=0\eqno(17)$$ for all
vector fields $F_1,F_2,...,F_{c+1}$ which are commutators of length
from $2$ to $\tau+1$ of vector fields in $\Gamma(M,N)$ and
$\Gamma(M,\overline N)$.  (Note: The lengths of $F_1,F_2,...,F_{c+1}$
may be different.)

(II) The dimension of $N(z,w)$ is $\ell-c$ for all $(z,w)\in M$ and
there exists a unique self-conjugate involutive $C^\infty$ bundle
${\Cal T}$ with fiber ${\Cal T}(z,w)$ for $(z,w)\in M$ such that
$N\oplus\overline N\subset{\Cal T}\subset\Bbb CTM$ and $\pi_{(z,w)}:
{\Cal T}(z,w)\rightarrow \Bbb CT_zS$ is an isomorphism for all
$(z,w)\in M$ (so the complex dimension of ${\Cal T}(z,w)$ is
$2\ell-c$).

(III) For every $(z^0,w^0)\in M$, there exists some neighborhood $U^M$
of $(z^0,w^0)$ in $M$ such that $U^M$ is foliated by graphs
of $C^\infty$ CR maps defined in a neighborhood $U$ of $z^0$ in $S$;
for such a map $f$, $f:U\rightarrow
\Bbb C^m$, and $q_k(z,f(z))=0$ for all $z\in U$ and all $k=1$ to $d$.
Furthermore, for $i=1$ to $d$ let $\phi_i$ be the $1$-form
$\partial_wq_i=\sum_{j=1}^m{\partial q_i\over
\partial w_j}dw_j$ and let $\phi$ be the $d$-form $\phi_1\wedge
\phi_2\wedge\phi_3\wedge...\wedge\phi_d$ which we write as
$\sum_{1\leq i_1<i_2<i_3<...<i_d\leq m}
\phi_{i_1,i_2,...,i_d}(z,w)dw_{i_1}\wedge dw_{i_2}\wedge
dw_{i_3}\wedge...\wedge dw_{i_d}$.  Then there exists a nonzero
$C^\infty $ function $C:U\rightarrow\Bbb C$ such that the $d$-form
$C(z)\phi(z,f(z))$ has coefficients which are CR functions on $U$.

(IV)  The dimension of $N(z,w)$ is $\ell-c$ for
all $(z,w)\in M$ and for every $k=1$ to $d$ we have that
$$\langle\partial p_1\wedge\partial p_2\wedge\partial p_3\wedge...\wedge
\partial p_c\wedge\partial q_k,F_1\wedge
F_2\wedge F_3\wedge...\wedge F_{c+1}\rangle=0\eqno(18)$$ for all
vector fields $F_1,F_2,...,F_{c+1}$ which are commutators of {\rm
arbitrary} length of vector fields in $\Gamma(M,N)$ and
$\Gamma(M,\overline N)$.  (Note: The lengths of $F_1,F_2,...,F_{c+1}$
may be different.)

If the above conditions hold, the graphs in (III) are integral
manifolds of $\Re\,{\Cal T}$ and the $(1,0)$ tangent space to the
graph of $f$ at $(z,f(z))$ is $N(z,f(z))$.  Suppose
$g:U\rightarrow\Bbb C^m$ is a CR map such that $q_k(z,g(z))=0$ for all
$z\in U$ and all $k=1$ to $d$, and such that $N$ (restricted to the
graph of $g$) is the $(1,0)$ tangent bundle to the graph of $g$.  Then
we must have $f=g$ for some such $f$ above.

\endproclaim


\noindent{\bf Note:}  Since from (2) the $\partial_w 
q_i(\cdot,\cdot)$ are linearly independent over $\Bbb C$ at every
point of $M$, we have that the form $\phi(z,w)$ is never zero.  If
$m=d$, $\phi(z,w)$ is of the form $\phi_{1,2,3,...,m}(z,w)dw_1\wedge
dw_2\wedge...\wedge dw_m$.  In this case, the property that
$\phi(z,w)$ satisfies is always automatic: let
$C(z)=1/\phi_{1,2,3,...,m}(z,f(z))$.  Thus in the case $d=m$ condition
(III) is merely a statement that $M$ is locally foliated by graphs of
CR mappings.

\noindent{\bf Proof:} It is obvious that (IV) implies
(I).  We shall prove $(I)\implies (II)\implies(III)\implies(IV)$.
Before proving any of these implications, we note that in each of them
we shall obtain the fact that for every $(z,w)\in M$ the complex
dimension of $N(z,w)$ is $\ell-c$, so
$\pi_{(z,w)}:N(z,w)\rightarrow\Nu^S_z$ is a (complex) linear
isomorphism.  We have this fact as an assumption in the first two
implications above, and in the third we shall prove it.  (We already
know this mapping is injective from Proposition 1; if the dimension of
$N(z,w)$ is $\ell-c$, then $\pi_{(z,w)}:N(z,w)\rightarrow\Nu^S_z$ is
also surjective since the complex dimension of $\Nu^S_z$ is $\ell-c$.)
In the proof of each implication, we shall make use of the fact that
the dimension of $N(z,w)$ is $\ell-c$ to construct a set of vector
fields as follows.  Take $(z,w)\in M$ and let ${\Cal T}(z,w)$ be the
complex vector space generated by the set of values $T(z,w)$, where
$T$ is a commutator of vector fields in $\Gamma(M,N)$ and
$\Gamma(M,\overline N)$ of length less than or equal to $\tau$.  We
have that ${\Cal T}(z,w)$ is the same as the complex vector space
generated by the set of values $T(z,w)$, where $T$ is instead a
commutator of vector fields in $\Gamma(U^M,N)$ and
$\Gamma(U^M,\overline N)$ of length less than or equal to $\tau$, and
$U^M$ is some open neighborhood of $(z,w)$ in $M$.  Given a point
$z_0\in S$, choose a neighborhood $G$ of $z^0$ in $S$ such that
there exist $s_1,s_2,...,s_{\ell-c}\in\Gamma(G,\Nu^S)$ such that
$\{s_j(z)\}_{j=1}^{\ell-c}$ forms a basis for $\Nu^S_z$ for $z\in G$.
Since $S$ is of finite type $\tau$ at $z^0$, we may choose vector
fields $t_1$, $t_2$,...,$t_c$ which are linear combinations of
commutators of the $s_i,\overline s_i$ such that the length of each
term of each $t_i$ is less than or equal to $\tau$ and such that
$$\{s_1,s_2,s_3,...,s_{\ell-c}\}\cup\{\overline s_1,\overline
s_2,\overline s_3,...,\overline s_{\ell-c}\}\cup
\{t_1,t_2,...,t_c\}\eqno(19)$$ constitutes a local basis for $\Bbb CTS$ near
$z^0$.  Furthermore, by appropriate change of basis we may assume that
each $t_i$ is a real tangent to $S$.  We assume by shrinking $G$ that
they constitute a local basis in $G$.  Let $G^M$ be as in (6) and let
$n_1,n_2,...,n_{\ell-c}$ be the $C^\infty$ vector fields in
$\Gamma(G^M,\Nu^M)$ determined by (13).  Then define $T_i$, $i=1$ to $c$
in the following manner: recall that $t_i$ is a linear combination of
Lie brackets of vector fields in the set $\{s_j\}
_{j=1}^{\ell-c}\cup\{\overline s_j\}_{j=1}^{\ell-c}$.  For every
appearance of an element $s_j$, $j=1$ to $\ell-c$ in that expression,
replace it with $n_j$ to form $T_i$.  (For example, if
$t_i=[s_1,\overline s_2]+[s_2,[s_3,\overline s_4]]+s_2$ then
$T_i=[n_1,\overline n_2]+[n_2,[n_3,\overline n_4]]+n_2$.)  Then we
have defined $T_i\in\Gamma(G^M,\Bbb CTM)$ and we claim that
$$\pi_{(z,w)}(T_i(z,w))=t_i(z).\eqno(20)$$ Suppose we select
$r_1,r_2,...$ in the set $\{s_j\} _{j=1}^{\ell-c}\cup\{\overline
s_j\}_{j=1}^{\ell-c}$ and $R_1,R_2,...$ are the corresponding elements
in $\Gamma(G^M,N)$ and $\Gamma(G^M,\overline N)$ such that
$\pi_{(z,w)}(R_j(z,w))=r_j(z)$ for all $(z,w)\in G^M$ and
$j=1,2,3,...$.  (See (13).)  To prove (20) it will suffice to prove
that
$$\pi_{(z,w)}( [R_j,[R_{j-1},[...[R_2,R_1]]..](z,w))=
[r_j,[r_{j-1},[...[r_2,r_1]]..](z)\eqno(21)$$ for $j\geq 1$ and
$(z,w)\in G^M$, since the $T_i$ are linear combinations of elements of
the form $[R_j,[R_{j-1},[...[R_2,R_1]]..]$.  We do this by induction
on $j$.  If $j=1$ then (21) follows from the definition of the $n_i$
in (13).  Now assume (21) is true for commutators
$R=[R_k,[R_{k-1},[...[R_2,R_1]]..]$ of length equal to $k$.  Then
$\pi_{(z,w)}(R(z,w))$ depends only on $z$ by the induction hypothesis,
so we let $r(z)=\pi_{(z,w)}(R(z,w))$.  Let $R'=[R_{k+1},R]$.  By the
induction hypothesis, the coefficients of $R_{k+1}$ and $R$ in
${\partial\over\partial z_\sigma}$, $\sigma=1$ to $\ell$ depend only
on $z\in\Bbb C^\ell$.  Thus when calculating the Lie bracket of
$R_{k+1}$ and $R$, the ${\partial\over\partial w_j}$ terms of each
annihilate the ${\partial\over\partial z_\sigma}$ coefficients of the
other.  Thus $\pi_{(z,w)}([R_{k+1},R](z,w))=[r_{k+1},r](z)$, as
desired, which implies that (21) holds for $j=k+1$.  By induction (21)
holds for all $j=1,2,3,...$, so (20) holds as well.  Since
$\pi_{(z,w)}(n_i(z,w))=s_i(z)$ and $\pi_{(z,w)}(\overline
n_i(z,w))=\overline s_i(z)$ for $i=1$ to $\ell-c$, and
$\pi_{(z,w)}(T_i(z,w))=t_i(z,w)$ for $i=1$ to $c$, we find that
$\pi_{(z,w)}:{\Cal T}(z,w)\rightarrow\Bbb CT_zS$ is surjective.
Furthermore, since (19) is a local basis for $\Bbb CTS$ in $G$, the
set of vectors $$\{n_i(z,w)\}_{i=1}^{\ell-c}\cup\{\overline
n_i(z,w)\}_{i=1}^{\ell-c}\cup\{T_i(z,w)\}_{i=1}^c\eqno(22)$$ is
linearly independent in $\Bbb CT_{(z,w)}M$ for all $(z,w)\in G^M$, as
the inverse image under $\pi_{(z,w)}$ of the basis
$\{s_i(z)\}_{i=1}^{\ell-c}\cup\{\overline
s_i(z)\}_{i=1}^{\ell-c}\cup\{t_i(z)\}_{i=1}^c$ for $\Bbb CT_zS$.

Now we assume (I).  There we simply assume that the dimension of
$N(z,w)$ is $\ell-c$, so the above vector fields
$\{n_i\}_{i=1}^{\ell-c}$ and $\{T_i\}_{i=1}^c$ can be constructed.  We
shall show that the quotient of ${\Cal T}(z,w)$ by
$N(z,w)\oplus\overline N(z,w)$ has complex dimension $c$, which we
recall is the CR codimension of $S$ in $\Bbb C^\ell$.  We will show
that the quotient ${\Cal T}(z,w)/(N(z,w)\oplus\overline N(z,w))$ is in
fact generated by the image in the quotient of $\{T_i(z,w)\}_{i=1}^c$;
this will show the quotient has dimension $c$, as desired.  The
(surjective) mapping $\pi_{(z,w)}:{\Cal T}(z,w)\rightarrow\Bbb CT_zS$
can be composed with the (surjective) quotient mapping from $\Bbb
CT_zS\rightarrow \Bbb CT_zS/(\Nu^S_z\oplus\overline \Nu^S_z)$ to form
another surjective mapping $\tilde\pi_{(z,w)}:{\Cal
T}(z,w)\rightarrow\Bbb CT_zS/(\Nu^S_z\oplus\overline \Nu^S_z)$.  Since
$\tilde\pi_{(z,w)}(n_i(z,w))=s_i(z)+\Nu^S_z\oplus\overline
\Nu^S_z =0+\Nu^S_z\oplus\overline \Nu^S_z$ and for similar
reasons $\tilde\pi_{(z,w)}(\overline n_i(z,w))$ is zero, we find that
$\tilde\pi_{(z,w)}$ factors through the quotient ${\Cal
T}(z,w)/(N(z,w) \oplus\overline N(z,w))$ to form a surjective mapping
from ${\Cal T}(z,w)/(N(z,w)\oplus\overline N(z,w))$ to $\Bbb CT_zS/
(\Nu^S_z\oplus\overline\Nu^S_z)$.  Since the complex dimension of the
latter space is $c$ and the mapping is surjective, the dimension of
${\Cal T}(z,w)/(N(z,w)\oplus\overline N(z,w))$ is greater than or
equal to $c$.  We proceed to show it does not exceed $c$.

\proclaim{Lemma 3} Assume that the conditions of Theorem 1 hold and that 
(I) holds.  Let $T(z,w)$ be a commutator of vector fields in the set
$\cup_{i=1}^{\ell-c}\{n_i,\overline n_i\}$.  We claim that there exist
$C^\infty $ functions $a_i:G\rightarrow\Bbb C$, $i=1$ to $c$,
such that
$$T+\sum_{i=1}^ca_iT_i\in\Gamma(G^M,N\oplus\overline
N).\eqno(23)$$

\endproclaim
\noindent (Note that by the definition of $G^M$ in (6), we may regard $a_i$
as a function on $G^M$ as well.)
\noindent{\bf Proof:} 

We shall first show that for some functions $a_i$
$$T+\sum_{i=1}^ca_iT_i\in\Gamma(G^M,\Nu^M\oplus\overline \Nu^M).
\eqno (24)$$ 
To do this, recall that the coefficients of $n_i$ in
${\partial\over\partial z_k}$ depend only on $z$ for all $i=1$ to
$\ell-c$ and $k=1$ to $\ell$, so from (21) the same may be said of the
coefficients of $T(z,w)$; thus $\pi_{(z,w)}(T(z,w))$ is a well-defined
function of $z$, say $s(z)$.  Because $s(z)\in\Bbb CT_zS$ and because
of the existence of the local basis (19) we may choose $C^\infty$
functions $a_i:G\rightarrow\Bbb C$, $i=1$ to $c$ such that
$s(z)+\sum_{i=1}^ca_i(z)t_i(z)\in
\Nu^S_z\oplus\overline\Nu^S_z$ for all $z\in G$.  
Then consider the vector field $T+\sum_{i=1}^ca_iT_i\in\Gamma(G^M,
\Bbb CTM)$.  We have that for all $j=1$ to $c$,
$$\eqalign{\langle\partial
p_j,T+\sum_{i=1}^ca_iT_i\rangle(z,w)&=\langle\partial
p_j(z),\pi_{(z,w)}(T(z,w)+\sum_{i=1}^ca_i(z)T_i(z,w))\rangle\cr
&=\langle\partial
p_j,s+\sum_{i=1}^ca_it_i\rangle(z)=0\cr}\eqno(25)$$ for all
$(z,w)\in G^M$ since $p_j$ depends only on $z$ and
$s(z)+\sum_{i=1}^ca_i(z)t_i(z)\in
\Nu^S_z\oplus\overline\Nu^S_z$ for all $z\in G$.

We claim that it is also then true that $$\langle\partial
q_k,T+\sum_{i=1}^ca_iT_i\rangle(z,w)=0\eqno(26)$$ for $k=1$ to $d$ and
$(z,w)\in G^M$.  We now use the definition of the action of 
a $c+1$-form on an element of the $c+1^{st}$ exterior algebra
of the tangent space (see [Bo, p. 10]) to calculate
$$\langle\partial p_1\wedge\partial p_2\wedge\partial
p_3\wedge...\wedge\partial p_c\wedge\partial q_k,T_1\wedge T_2\wedge
T_3\wedge...\wedge T_c\wedge (T+\sum_{i=1}^ca_iT_i)\rangle(z,w)$$
for $(z,w)\in G^M$.  Any term involving $\langle\partial
p_i,T+\sum_{i=1}^ca_iT_i\rangle$ is zero by (25), so we have from
(17) that
$$\eqalign{0&=\langle\partial p_1\wedge\partial p_2\wedge\partial
p_3\wedge...\wedge\partial p_c\wedge\partial q_k,T_1\wedge T_2\wedge
T_3\wedge...\wedge T_c\wedge (T+\sum_{i=1}^ca_iT_i)\rangle(z,w)\cr
&=\langle\partial p_1\wedge\partial p_2\wedge\partial
p_3\wedge...\wedge\partial p_c,T_1\wedge T_2\wedge T_3\wedge...\wedge
T_c\rangle\langle
\partial q_k,T+\sum_{i=1}^ca_iT_i\rangle(z,w)\cr}$$
for $(z,w)\in G^M$.  We claim that $\langle\partial
p_1\wedge\partial p_2\wedge\partial p_3\wedge...\wedge\partial
p_c,T_1\wedge T_2\wedge T_3\wedge...\wedge T_c\rangle(z,w)\neq 0$ for
$(z,w)\in G^M$; this will show that (26) holds, as desired.  To see
why the claim holds, note that $\langle\partial p_1\wedge\partial
p_2\wedge\partial p_3\wedge...\wedge\partial p_c,T_1\wedge T_2\wedge
T_3\wedge...\wedge T_c\rangle(z,w)$ is equal to the determinant
of the $c\times c$ matrix whose $i,j$ component is $\langle\partial
p_i,T_j\rangle(z,w)$.  If at some $(z^0,w^0)\in G^M$ this determinant is
zero then the columns are linearly dependent over $\Bbb C$ , so for
some $\zeta_j\in \Bbb C$ (which are not all zero) and all $i=1$ to
$c$, $\sum_{j=1}^c\zeta_j\langle\partial p_i, T_j\rangle(z^0,w^0)=0$, so
$\langle\partial p_i,\sum_{j=1}^c\zeta_j T_j\rangle(z^0,w^0)=0$ and
$\langle\partial p_i,\sum_{j=1}^c\zeta_j t_j\rangle(z^0)=0$ for $i=1$ to
$c$.  By definition of $\Nu^S_{z^0}$ we have $\sum_{j=1}^c\zeta_j
t_j(z^0)\in\Nu^S_{z^0}\oplus\overline\Nu^S_{z^0}$, so for some $a_j,b_j\in \Bbb
C$, $j=1$ to $\ell-c$, we obtain $\sum_{j=1}^c\zeta_j
t_j(z^0)=\sum_{j=1}^{\ell-c}a_j s_j(z^0)+\sum_{j=1}^{\ell-c}b_j \overline
s_j(z^0)$.  Since the set $\{s_i(z^0)\}_{i=1}^{\ell-c}\cup\{\overline
s_i(z^0)\}_{i=1}^{\ell-c}\cup \{t_i(z^0)\}_{i=1}^{c}$ is linearly independent,
we find that $\zeta_j=0$ for $j=1$ to $c$.  This contradicts the
definition of the $\zeta_j$, so the claim holds.

By (25) and (26) we conclude that (24) holds.  We proceed to show that
(23) also holds.  To see this, we observe that there exists
$h\in\Gamma(G^M,\Nu^M)$ such that $T+\sum_{i=1}^ca_iT_i+\overline
h\in\Gamma(G^M,\Nu^M)$.  Then for $v\in\Gamma(G^M,V)$,
$[T+\sum_{i=1}^ca_iT_i+\overline h,\overline
v]=[T+\sum_{i=1}^ca_iT_i,
\overline v]+[\overline h,\overline v]$.  The first of these
last two terms belongs to $\Gamma(G^M,\Nu^M\oplus\overline \Nu^M)$ by
Lemma 1 and the second belongs to $\Gamma(G^M,\Nu^M\oplus\overline
\Nu^M)$ since $\overline\Nu^M$ is involutive.

Recalling that every element of $V(z,w)$ is the value at $(z,w)$ of
some such $v\in\Gamma(M,V)$ and recalling the definition of $N$, we
may write ${\Cal
L}^M_{(z,w)}(T(z,w)+\sum_{i=1}^ca_i(z)T_i(z,w)+\overline
h(z,w),v_{zw})=0$ for all $v_{zw}\in V(z,w)$ and all $(z,w)\in M$.
Thus we may write that $T+\sum_{i=1}^ca_i(z)T_i+\overline
h=n'\in\Gamma(G^M,N).$ If we solve for $h$ here we can use a similar
argument to show that $h\in\Gamma(G^M,N)$, so (23)
holds and Lemma 3 is proven.$\hfill\square$

Lemma 3 shows that the quotient of ${\Cal T}(z,w)$ by
$N(z,w)\oplus\overline N(z,w)$ has complex dimension less than or
equal to $c$.  We already know the dimension is greater than or equal
to $c$, so it is exactly $c$.  The set in (22) is linearly
independent for $(z,w)\in G^M$ as observed earlier and Lemma 3 shows
that it spans ${\Cal T}(z,w)$ for $(z,w)\in G^M$.  Since the $T_i$ are
chosen smoothly in the neighborhood $G$, the set of vector fields
$$\{n_i\}_{i=1}^{\ell-c}\cup\{\overline
n_i\}_{i=1}^{\ell-c}\cup\{T_i\}_{i=1}^c$$ is a local basis for ${\Cal
T}$ over $G^M$.  Since sets $G^M$ cover $M$, the ${\Cal T}(z,w)$ form a
$C^\infty$ bundle on $M$ of complex dimension $2\ell-c$ which we call
${\Cal T}$.

Next we show that ${\Cal T}$ is involutive.  Suppose we select an
arbitrary $(z^0,w^0)\in M$, an open neighborhood $U^M$ of $(z^0,w^0)$
in $M$ and sections $R_1,R_2,...R_{2\ell-c}$ of $\Gamma(U^M,{\Cal T})$
such that $\{R_i(z,w)\}_{i=1}^{2\ell-c}$ is a basis for ${\Cal
T}(z,w)$ for $(z,w)\in U^M$. Then to show that ${\Cal T}$ is
involutive it suffices to show that Lie brackets of the $R_i$ belong
to $\Gamma(U^M,{\Cal T})$.  In fact, we can let $U^M$ be the $G^M$ be
defined in (6) and let the set $\{R_i\}_{i=1}^{2\ell-c}$ be
$\{n_i\}_{i=1}^{\ell-c}\cup\{\overline
n_i\}_{i=1}^{\ell-c}\cup\{T_k\}_{k=1}^c$ in $G^M$.  It will be enough
to show that for $(z,w)\in G^M$, we have that $[n_i, n_j] (z, w)$,
$[n_i, \overline n_j] (z, w)$, $[T_k,n_i] (z, w)$, and $[T_k,\overline
n_i] (z, w)$ all belong to ${\Cal T}(z,w)$.  These are all
consequences of Lemma 3.  Thus ${\Cal T}$ is involutive.  We have that
${\Cal T}$ is self-conjugate because the local basis
$\{n_i(z,w)\}_{i=1}^{\ell-c}\cup\{\overline
n_i(z,w)\}_{i=1}^{\ell-c}\cup
\{T_i(z,w)\}_{i=1}^{c}$ for ${\Cal T}(z,w)$ is self-conjugate.

We need to show that $\pi_{(z,w)}:{\Cal T}(z,w)\rightarrow\Bbb
CT_zS$ is an isomorphism for every $(z,w)\in M$.  Once again it
will suffice to assume that $(z,w)\in G^M$ where $G,G^M$ are as
before.  We have that $\pi_{(z,w)}(n_i(z,w))=s_i(z)$ ($i=1$ to
$\ell-c$), $\pi_{(z,w)}(\overline n_i(z,w))=\overline s_i(z)$ ($i=1$
to $\ell-c$), and $\pi_{(z,w)}(T_i(z,w))=t_i(z)$ ($i=1$ to $c$).
Since $\{n_i(z,w)\}_{i=1}^{\ell-c}\cup\{\overline
n_i(z,w)\}_{i=1}^{\ell-c}\cup
\{T_i(z,w)\}_{i=1}^{c}$ and $\{s_i(z)\}_{i=1}^{\ell-c}\cup\{\overline
s_i(z)\}_{i=1}^{\ell-c}\cup
\{t_i(z)\}_{i=1}^{c}$ are bases, respectively, for ${\Cal T}(z,w)$
and $\Bbb CT_zS$, the map $\pi_{(z,w)}:{\Cal
T}(z,w)\rightarrow\Bbb CT_zS$ is indeed an isomorphism.

Now we prove the uniqueness statement of $(II)$.  Every involutive
bundle ${\Cal T}'$ on $G^M$ between $N\oplus\overline N$ and $\Bbb
CTM$ must contain $T_i(z,w)$ in its fiber over $(z,w)$ (for $i=1$ to
$c$) since the $T_i$ are linear combinations of commutators of the
$n_j,\overline n_j$.  Thus ${\Cal T}'(z,w)$ must contain ${\Cal
T}(z,w)$.  For $\pi_{(z,w)}:{\Cal T}'(z,w)\rightarrow\Bbb
CT_{(z,w)}S$ to be an isomorphism, ${\Cal T}'(z,w)$ must have have
complex dimension $2\ell-c$, so it must be no bigger than ${\Cal T}$
over $G^M$, since the complex dimension of ${\Cal T}$ is $2\ell-c$.
This shows that ${\Cal T}'(z,w)={\Cal T}(z,w)$ for $(z,w)\in G^M$.
Since open sets such as $G^M$ cover $M$, we must have ${\Cal
T}'(z,w)={\Cal T}(z,w)$ for $(z,w)\in M$; this proves the uniqueness
statement of $(II)$ and concludes the proof that (I) implies (II).

Now we show that $(II)\implies (III)$.  We have that $\pi_{(z,w)}:
N(z,w)\rightarrow \Nu^S_z$ is an isomorphism and the complex dimension
of $N(z,w)$ is $\ell-c$.  Fix $(z^0,w^0)\in M$ and let
$\{s_i\}_{i=1}^{\ell-c},\{n_i\}_{i=1}^{\ell-c}$ and sets $G,G^M$ be as
defined in (5),(6),(13).  Because the dimension of $N(z,w)$ is
$\ell-c$, we may define $\{t_i\}$ as in (19) and $\{T_i\}$ as in
(20).  Because $\pi_{(z,w)}:{\Cal T}(z,w)\rightarrow
\Bbb CT_zS$ is an isomorphism and $\{s_i(z)\}_{i=1}^{\ell-c}\cup\{\overline
s_i(z)\}_{i=1}^{\ell-c}\cup
\{t_i(z)\}_{i=1}^{c}$ is a basis for $\Bbb CT_zS$, the set 
$\{n_i(z,w)\}_{i=1}^{\ell-c}\cup\{\overline n_i(z,w)\}_{i=1}^{\ell-c}\cup
\{T_i(z,w)\}_{i=1}^{c}$ is a basis for ${\Cal T}(z,w)$ for $(z,w)\in G^M$
and $G$ sufficiently small.  Let us define $\Re\,{\Cal T}(z,w)$ to be
the vector space $\{R_{zw}+\overline {R_{zw}}|R_{zw}\in {\Cal
T}(z,w)\}$.  We check that $\Re\,{\Cal T}$ constitutes a real
involutive vector bundle over $M$ of real dimension $2\ell-c$.  It
will suffice to check that $\Re\,{\Cal T}$ constitutes a real
involutive vector bundle over $G^M$ of real dimension $2\ell-c$ for
every $G^M$ defined above, since such $G^M$ cover $M$.  Since
$\{n_j(z,w):j=1,2,...,\ell-c\}\cup\{\overline
n_j(z,w):j=1,2,...,\ell-c\}
\cup\{T_j(z,w):j=1,2,...,c\}$
is a complex basis for ${\Cal T}(z,w)$, for $(z,w)\in G^M$ (recalling
that the $T_j$ are real vector fields), $\{n_j(z,w)+\overline
n_j(z,w):j=1,2,...,\ell-c\}\cup
\{in_j(z,w)+\overline {in_j(z,w)}|j=1,2,...,\ell-c\}\cup
\{T_j(z,w):j=1,2,...,c\}$ is a real basis for $\Re\,{\Cal T}(z,w)$,
so the (real) dimension of $\Re\,{\Cal T}(z,w)$ is $2\ell-c$ and
$\Re\,{\Cal T}$ is a bundle.  Also, $\Re\,{\Cal T}$ is involutive as
the real part of an involutive bundle.  By the Frobenius theorem (see
[Wa]) $M$ is foliated near $(z^0,w^0)\in M$ by $C^\infty$ integral
manifolds for $\Re\,{\Cal T}$.  The complexified tangent space to such
a manifold at $(z,w)\in M$ must equal ${\Cal T}(z,w)$.

We know that $\pi_{(z,w)}$ is injective on ${\Cal T}(z,w)$; hence it
is injective on $\Re\,{\Cal T}(z,w)$ also.  By the inverse function
theorem, near $(z^0,w^0)$ the integral manifolds of $\Re\,{\Cal T}$
are graphs over a fixed open subset of $S$, so we write that a
neighborhood of $(z^0,w^0)$ is foliated by graphs of mappings on an
open subset $U$ of $S$ where $z^0\in U$.  The complexified tangent
bundles to these manifolds must equal ${\Cal T}$.  Let $f$ be the
function on an open subset $U$ of $S$ such that $z^0\in U$ and such
that the graph of $f$ is the integral manifold of $\Re\,{\Cal T}$
through $(z^0,w^0)$.  Use $M^f$ to denote the graph of $f$.  The
$(1,0)$ tangent space to $M^f$ at $(z,f(z))$ includes $N(z,f(z))$
(since the whole complexified tangent space to $M^f$ at $(z,f(z))$ is
${\Cal T}(z,f(z))$.)  Any $(1,0)$ tangent to $M^f$ at $(z,f(z))$
projects to $\Nu^S_z$; since projection is injective on ${\Cal
T}(z,f(z))=\Bbb CT_{(z,f(z))}M^f$, that tangent must belong to
$N(z,f(z))$.  This shows that the $(1,0)$ tangent space to $M^f$ at
$(z,f(z))$ is $N(z,f(z))$.

Then the differential of the mapping $F:U\rightarrow M^f$ such that
$F(z)=(z,f(z))$ from $U$ to $M^f$ carries $\Nu^S$ to $N$: for $s'_z\in
\Nu^S_z$, if $d_zF(s'_z)$ is not in $N(z,f(z))$, then
$\pi_{(z,f(z))}(d_zF(s'_z))=s'_z$ is not in $\Nu^S_z$ (contradiction),
since $\pi_{(z,f(z))}$ maps $N(z,f(z))$ onto $\Nu^S_z$ and is injective
on ${\Cal T}(z,f(z))$.  Thus $F$ must be a CR mapping on $U$, so $f$
is as well.

We now proceed to prove the property that (III) asserts for the
$d$-form $\phi$.  Let $s\in\Gamma(U,\Nu^S)$.  If $K$ is a function on
$S$ then we write $s_z\{K\}$ to denote the action of $s$ on $K$ at
$z\in S$.  (Then we regard $s_z$ as an element of $\Nu^S_z$ and
identify $s(z)$ with $s_z$.  We use $\overline s_z$ similarly.)  For
convenience, we note that, since $f$ is CR, then $$n(z)\equiv
s(z)+\sum_{j=1}^m s_z\{f_j\}{\partial
\over \partial w_j}\eqno(27)$$ belongs to $N(z,f(z))$, as the element
of the $(1,0)$ tangent space $N(z,f(z))$ to the graph of $f$ at
$(z,f(z))$ which projects to $s(z)$.

We note that the function $\phi_{i_1,i_2,...,i_d}(z,w)$ (defined in
property $(III)$) is the determinant of the $d\times d$ matrix with
$(j,k)$ entry ${\partial q_j\over\partial w_{i_k}}(z,w)$.  Since at
every $(z,w)\in M$, the set $\{\partial_wq_i(z,w):i=1,2,...,d\}$ is
linearly independent (from (2)), $\phi(z,w)$ is nonzero for $(z,w)\in
M$.  Thus for some integers $i_1,i_2,...,i_d$, $1\leq
i_1<i_2<...<i_d\leq m$, we have that
$\phi_{i_1,i_2,...,i_d}(z^0,w^0)\neq 0$.  We claim that it suffices to
take $C(z)=1/\phi_{i_1,i_2,...,i_d}(z,f(z))$ for $z$ in a possibly
shrunken neighborhood $U$ of $z^0\in S$.  Assume without loss of
generality that $C$ is defined on all of $U$ (by shrinking $U$).  Thus
we will show that for $\overline s \in
\Gamma(U,\overline \Nu^S)$ and integers $j_k$, $1\leq j_1<j_2<
j_3<...<j_d\leq
m$,

$$\overline s_z\left\{{\left |\matrix {\partial q_1\over\partial
w_{j_1}} (\cdot,f(\cdot)) & {\partial q_1\over\partial
w_{j_2}}(\cdot,f(\cdot)) & \dots &{\partial q_1\over\partial
w_{j_d}}(\cdot,f(\cdot))\\ {\partial q_2\over\partial
w_{j_1}}(\cdot,f(\cdot)) & {\partial q_2\over\partial
w_{j_2}}(\cdot,f(\cdot)) & \dots &{\partial q_2\over\partial
w_{j_d}}(\cdot,f(\cdot))\\
\vdots&&\dots&\vdots\\
{\partial q_d\over\partial w_{j_1}}(\cdot,f(\cdot)) & {\partial
q_d\over\partial w_{j_2}}(\cdot,f(\cdot)) & \dots &{\partial
q_d\over\partial w_{j_d}}(\cdot,f(\cdot))
\endmatrix \right |
\over
\left |\matrix {\partial q_1\over\partial w_{i_1}}(\cdot,f(\cdot)) &
{\partial q_1\over\partial w_{i_2}}(\cdot,f(\cdot)) & \dots &{\partial
q_1\over\partial w_{i_d}}(\cdot,f(\cdot))\\ {\partial q_2\over\partial
w_{i_1}}(\cdot,f(\cdot)) & {\partial q_2\over\partial
w_{i_2}}(\cdot,f(\cdot)) & \dots &{\partial q_2\over\partial
w_{i_d}}(\cdot,f(\cdot))\\
\vdots&&\dots&\vdots\\
{\partial q_d\over\partial w_{i_1}}(\cdot,f(\cdot)) & {\partial
q_d\over\partial w_{i_2}}(\cdot,f(\cdot)) & \dots &{\partial
q_d\over\partial w_{i_d}}(\cdot,f(\cdot))
\endmatrix \right |
}\right\}\eqno(28)$$ is identically zero for $z\in U$.  Let $A(z)$ be
the matrix with $(\sigma,k)$ entry ${\partial q_\sigma\over\partial w_
{j_k}}(z,f(z))$, let $B(z)$ be the matrix with $(\sigma,k)$ entry
${\partial q_\sigma\over\partial w_{i_k}}(z,f(z))$ and let $\det A(z),\det
B(z)$ be their determinants.  Then (28) is equal to
$${\det B(z)\overline s_z\{\det A\}- \det A(z)\overline s_z\{\det B\} 
\over [\det(B(z))]^2}.\eqno(29)$$
Noting that $\phi_{j_1,j_2,j_3,...,j_d}(z,f(z))=\det A(z)$ and
$\phi_{i_1,i_2,i_3,...,i_d}(z,f(z))=\det B(z)$ then what we have to
prove is the following lemma, which is stated as a lemma because we
will need it later.

\proclaim{Lemma 4}
Suppose that the conditions of Theorem 1 are satisfied, (II) holds,
$U\subset S$, $f:U\rightarrow\Bbb C^m$ is a CR map whose graph is an
integral manifold for $\Re\,{\Cal T}$ and $s\in\Gamma(U,\Nu^S)$.  Then
for all $z\in U$ and any $d$-tuples
$\{i_k\}_{k=1}^d$,$\{j_k\}_{k=1}^d$ such that $1\leq
i_1<i_2<...<i_d\leq m$ and $1\leq j_1<j_2<...<j_d\leq m$ we have
$$\phi_{i_1,i_2,i_3,...,i_d}(z,f(z))\overline
s_z\{\phi_{j_1,j_2,j_3,...,j_d}(\cdot,f(\cdot))\}-\phi_{j_1,j_2,j_3,...,j_d}
(z,f(z))\overline
s_z\{\phi_{i_1,i_2,i_3,...,i_d}(\cdot,f(\cdot))\}=0\eqno(30)$$ where
$\phi_{i_1,i_2,i_3,...,i_d},\phi_{j_1,j_2,j_3,...,j_d}$ are the
functions defined in (III).
\endproclaim

\noindent{\bf Proof:}
We let $A,B$ be the matrices defined above.  Then we must show that
$$\det B(z)\overline s_z\{\det A\}- \det A(z)\overline s_z\{\det
B\}=0.\eqno(31)$$ We have $$\overline s_z\{\det A\}=\sum_{i=1}^d
\left |\matrix A_{11}(z)& A_{12}(z) & \dots & A_{1d}(z)\\ A_{21}(z)&
A_{22}(z)& \dots & A_{2d}(z)\\
\vdots&&\dots&\vdots\\ \overline s_z\{A_{i1}\}&\overline
s_z\{A_{i2}\}&\dots&\overline s_z\{A_{id}\}\\
\vdots&&\dots&\vdots\\ A_{d1}(z) & A_{d2}(z) & \dots &A_{dd}(z)
\endmatrix\right|\eqno(32)$$
where we calculate $$\overline s_z\{A_{i,k}\}=
\sum_{\sigma=1}^\ell{\partial q_i\over \partial \overline z_\sigma\partial 
w_{j_k}}(z,f(z))\overline
s_z\{\overline z_\sigma\}+\sum_{\sigma=1}^m{\partial q_i\over 
\partial \overline 
w_\sigma\partial w_{j_k}}(z,f(z))\overline s_z\{\overline f_\sigma\}.$$
We observe that if we let 
$$v^\sigma(z)=\det\left(\matrix A_{11}(z)& A_{12}(z) & \dots & A_{1d}(z)\\
A_{21}(z)& A_{22}(z)& \dots & A_{2d}(z)\\ \vdots&&\dots&\vdots\\
{\partial\over\partial w_{j_1}}&{\partial\over\partial
w_{j_2}}&\dots&{\partial\over\partial w_{j_d}}\\
\vdots&&\dots&\vdots\\ A_{d1}(z) & A_{d2}(z) & \dots &A_{dd}(z)
\endmatrix\right),$$
(where the row with the ${\partial\over\partial w_{j_k}}$, ${k=1}$ to $d$
is the $\sigma^{th}$ row) then for $i=1$ to $d$, $i\neq \sigma$,
$$\langle\partial q_i(z,f(z)),v^\sigma(z)\rangle= \det\left(\matrix A_{11}(z)&
A_{12}(z) & \dots & A_{1d}(z)\\ A_{21}(z)& A_{22}(z)& \dots & A_{2d}(z)\\
\vdots&&\dots&\vdots\\ {\partial q_i\over\partial w_{j_1}}(z,f(z))&
{\partial
q_i\over\partial w_{j_2}}(z,f(z))&\dots&{\partial
q_i\over\partial w_{j_d}}(z,f(z))\\
\vdots&&\dots&\vdots\\ A_{d1} (z)& A_{d2} (z)& \dots &A_{dd}(z)
\endmatrix\right)=0\eqno(33)$$
for all $z\in U$ because the $i^{th}$ and $\sigma^{th}$ rows of the
determinant are the same.  For $i=\sigma$, we have $$\langle\partial
q_i(z,f(z)),v^\sigma(z)\rangle=\det A(z).\eqno(34)$$ Write
$v^\sigma=\sum_{k=1}^d v^\sigma_k{\partial\over\partial w_{j_k}}$.
Combining these observations, we find that (32) guarantees that
$$\eqalign{\overline s_z\{\det A\}&=\sum_{i=1}^d\sum_{k=1}^d\overline
s_z\{A_{i,k}\}v^i_k\cr &=\sum_{i=1}^d\sum_{k=1}^d
\left (\sum_{\sigma=1}^\ell{\partial q_i\over \partial \overline
z_\sigma\partial w_{j_k}}(z,f(z))\overline s_z\{\overline
z_\sigma\}+\sum_{\sigma=1}^m{\partial q_i\over
\partial \overline 
w_\sigma\partial w_{j_k}}(z,f(z))\overline s_z\{\overline f_\sigma\}\right )
v^i_k\cr &=\sum_{i=1}^d\left(\sum_{k=1}^d
\sum_{\sigma=1}^\ell{\partial q_i\over \partial \overline
z_\sigma\partial w_{j_k}}(z,f(z))\overline s_z\{\overline
z_\sigma\}v^i_k+\sum_{k=1}^d\sum_{\sigma=1}^m{\partial q_i\over
\partial \overline 
w_\sigma\partial w_{j_k}}(z,f(z))\overline s_z\{\overline f_\sigma\}
v^i_k\right)\cr   &=\sum_{i=1}^d\langle\overline\partial\partial
q_i(z,f(z)),
\overline {n(z)}\wedge v^i(z)\rangle,\cr}\eqno(35)$$ where $n(z)$ was
defined in (27).  A similar argument
shows that
$$\overline s_z\{\det
B\}=\sum_{i=1}^d\langle\overline\partial\partial q_i(z,f(z)),
\overline {n(z)}\wedge u^i(z)\rangle.\eqno(36)$$ where
$$u^\sigma(z)=\det\left(\matrix B_{11}(z)& B_{12}(z) & \dots & B_{1d}(z)\\
B_{21}(z)& B_{22}(z)& \dots & B_{2d}(z)\\ \vdots&&\dots&\vdots\\
{\partial\over\partial w_{i_1}}&{\partial\over\partial
w_{i_2}}&\dots&{\partial\over\partial w_{i_d}}\\
\vdots&&\dots&\vdots\\ B_{d1}(z) & B_{d2}(z) & \dots &B_{dd}(z)
\endmatrix\right),$$
and where once again, the row with the ${\partial\over\partial w_{i_k}}$
is the $\sigma^{th}$ row.  As before, $$\langle\partial
q_j(z,f(z)),u^\sigma(z)\rangle=0\eqno(37)$$ if $j\neq \sigma$ and
$$\langle\partial q_j(z,f(z)),u^\sigma(z)\rangle=\det B(z)\eqno(38)$$ if
$j=\sigma$.  Then (31) is equal to
$$\sum_{i=1}^d\langle\overline\partial\partial q_i(z,f(z)), \overline
{n(z)}\wedge((\det B(z))v^i(z)-(\det
A(z))u^i(z))\rangle,\eqno(39)$$ where we note that for $j=1$ to $d$,
$$\eqalign{&\langle\partial q_j(z,f(z)),(\det B(z))v^i(z)-(\det
A(z))u^i(z)\rangle\cr=&
\det B(z)\langle\partial
q_j(z,f(z)),v^i(z)\rangle-\det A(z)\langle\partial q_j(z,f(z)),
u^i(z)\rangle\cr=&0\cr}$$ from (33),(34),(37) and (38).  By definition
of $V(z,f(z))$, $(\det B(z))v^i(z)-(\det A(z))u^i(z)\in V(z,f(z))$;
this implies that (39) is identically zero for $z\in U$, since
$n(z)\in N(z,f(z))$.  Thus (31) holds, so (30) holds as well, as
desired.$\hfill\square$

Now that we know Lemma 4 holds, we have that that (29) and (28) are
also zero for all $z\in U$, as desired.  This proves that last
component of (III), and hence we now know that (II) implies (III).

Now we assume (III) and prove (IV).  Then we know that there exists a
nonzero complex-valued function $C(z)$ defined on $U$ such that
$$C(z)\phi(z,f(z))=\sum_{1\leq i_1<i_2<i_3<...<i_d\leq
m}C(z)\phi_{i_1,i_2,...,i_d}(z,f(z)) dw_{i_1}\wedge
dw_{i_2}\wedge...\wedge dw_{i_d}$$ is a $d$-form with coefficients
which are CR functions for $z\in U$.  If $s\in\Gamma(U,\Nu^S)$ then
$\overline s_z\{C(\cdot)\phi(\cdot,f(\cdot))\}=0$ for all $z\in U$.
Then, using the product rule for the differentiation of wedge
products, $$\eqalign{0=&\overline s_z\{C\}\phi(z,f(z))\cr
+&C(z)\overline s_z
\{(\sum_{j=1}^m{\partial q_1\over\partial
w_j}(\cdot,f(\cdot))dw_j)\wedge(\sum_{j=1}^m{\partial q_2\over\partial
w_j}(\cdot,f(\cdot))dw_j)\wedge...\wedge(\sum_{j=1}^m{\partial q_d\over\partial
w_j}(\cdot,f(\cdot))dw_j)\}\cr =&\overline
s_z\{C\}\phi(z,f(z))\cr +&C(z)\sum_{i=1}^d
\biggl (
(\sum_{j=1}^m{\partial q_1\over\partial
w_j}(z,f(z))dw_j)\wedge(\sum_{j=1}^m{\partial q_2\over\partial
w_j}(z,f(z))dw_j)\wedge...\wedge(\sum_{j=1}^m{\partial
q_{i-1}\over\partial w_j}(z,f(z))dw_j)\wedge\cr&\hskip 1cm
(\sum_{j=1}^m\sum_{\sigma=1}^\ell{\partial^2q_i\over\partial\overline
z_\sigma\partial w_j}(z,f(z))\overline s_z\{\overline
z_\sigma\}dw_j+\sum_{j,\sigma=1}^m{\partial^2q_i\over
\partial\overline w_\sigma\partial w_j}(z,f(z))\overline s_z \{\overline 
f_\sigma\} dw_j)\cr & \hskip 1cm \wedge(\sum_{j=1}^m{\partial
q_{i+1}\over\partial w_j}(z,f(z))dw_j)
\wedge...\wedge(\sum_{j=1}^m{\partial q_d\over\partial
w_j}(z,f(z))dw_j) 
\biggr )
.\cr}\eqno(40)$$

For fixed $z\in U$, we claim that there exist $v^1,v^2,...,v^d$ in
$\Bbb CT_{f(z)}(M_z)$ such that for $i,j=1$ to $d$,
$$\langle\partial q_i(z,f(z)),v^j\rangle=\delta_{ij},\eqno(41)$$ where
$\delta_{ij}=0$ if $i\neq j$ and $\delta_{ij}=1$ if $i=j$.  Recall
from (2) that $M_z$ is a generic CR manifold in $\Bbb C^m$ of CR
codimension $d$.  The $1$-form $\partial q_i(z,f(z))$ induces a linear
functional on $\Bbb CT_{f(z)}(M_z)$ in the following manner: for
$x\in\Bbb CT_{f(z)}(M_z)$, map $x\mapsto L_i(x)\equiv\langle\partial
q_i(z,f(z)),x
\rangle$.  Note that $V(z,f(z))=\Bbb CT^{1,0}_{f(z)}(M_z)$.
Since for all $i=1$ to $d$, $L_i(x)=0$ for $x\in
V(z,f(z))\oplus\overline{V(z,f(z))}\subset\Bbb CT_{f(z)}(M_z)$ the
mapping $L_i$ factors through the quotient $\Bbb CT_{f(z)}(M_z)/
(V(z,f(z))\oplus\overline{V(z,f(z))})$ to produce a linear functional
$\tilde L_i$.  Now $\Bbb CT_{f(z)}(M_z)/
(V(z,f(z))\oplus\overline{V(z,f(z))})$ is a complex vector space of
dimension $d$ since $M_z$ is generic.  We claim that the induced
linear functionals $\tilde L_i$ are a basis of the dual space of $\Bbb
CT_{f(z)}(M_z)/ (V(z,f(z))\oplus\overline{V(z,f(z))})$.  If there
exist $\zeta_i\in
\Bbb C$ (not all zero) such that $\sum_{i=1}^d\zeta_i
\tilde L_i$ is identically zero on $\Bbb CT_{f(z)}(M_z)/
(V(z,f(z))\oplus\overline{V(z,f(z))})$ then $\sum_{i=1}^d\zeta_i L_i$
is identically zero as a linear functional on $\Bbb CT_{f(z)}(M_z)$
i.e. $\langle
\sum_{i=1}^d\zeta_i\partial q_i(z,f(z)),x\rangle=0$ for all $x\in 
\Bbb CT_{f(z)}(M_z)$.  Let $\i$ be the imaginary unit and let 
$J:T_{f(z)}\Bbb C^m\rightarrow T_{f(z)}\Bbb C^m$ be
the complex structure mapping on $T_{f(z)}\Bbb C^m$.  Then, extending
$J$ to $\Bbb CT_{f(z)}(\Bbb C^m)$, we recall that $J({\partial\over\partial
w_j}) =\i{\partial\over\partial w_j}$ and $J({\partial\over\partial
\overline w_j}) =-\i{\partial\over\partial \overline w_j}$ for all $j$.  
We find that $\langle
\sum_{j=1}^d\zeta_j\partial q_j(z,f(z)),J(x)\rangle=\i\langle
\sum_{j=1}^d\zeta_j\partial q_j(z,f(z)),x\rangle=0$ for all $x\in 
\Bbb CT_{f(z)}(M_z)$ because $\sum_{j=1}^d\zeta_j\partial q_j(z,f(z))$ is a
$(1,0)$ form.  Thus $\sum_{j=1}^d\zeta_j \partial q_j(z,f(z))$ is zero
as a linear functional on $J(\Bbb CT_{f(z)}(M_z))$.  Since $M_z$ is
generic, $\Bbb CT_{f(z)}(M_z)+J(\Bbb CT_{f(z)}(M_z))= \Bbb
CT_{f(z)}\Bbb C^m$ (see [BER, p. 14]), so $\sum_{j=1}^d\zeta_j\partial
q_j(z,f(z))$ is zero as a linear functional on $\Bbb CT_{f(z)}\Bbb
C^m$: $\sum_{j=1}^d\zeta_j\partial q_j(z,f(z))=0$.  This is impossible
because $M_z$ is generic.  Thus the $d$ functionals $\{\tilde L_i\}$
above are linearly independent in the dual space of $\Bbb
CT_{f(z)}(M_z)/ (V(z,f(z))\oplus\overline{V(z,f(z))})$, so are a basis
of that dual space (which has dimension $d$).  Let $\{\tilde
v^i\}_{i=1}^d$ be a basis of $\Bbb CT_{f(z)}(M_z)/
(V(z,f(z))\oplus\overline{V(z,f(z))})$ dual to $\{\tilde L_i\}$ and
for every $i=1$ to $d$ let $v^i$ be an element of $\Bbb
CT_{f(z)}(M_z)$ whose image in the quotient $\Bbb CT_{f(z)}(M_z)/
(V(z,f(z))\oplus\overline{V(z,f(z))})$ is $\tilde v^i$.  The $\{v^i\}$
satisfy (41), as desired.  Note that this implies that
$\langle\partial_w q_i(z,f(z)),v^j\rangle=\delta_{ij}$ since $v^j$ has
no terms involving ${\partial\over\partial z_k}$ for any $k$.

Let $v$ be an arbitrary element of $V(z,f(z))$.  If $R$ is the
rightmost side of (40) then we calculate $\langle R,v\wedge v^1\wedge
v^2\wedge...\wedge\widehat v^k\wedge ...\wedge v^d\rangle$, where
$\widehat v^k$ indicates $v^k$ is not in the wedge product.  Since
$\langle\partial_w q_j(z,f(z)),v\rangle=\langle\partial
q_j(z,f(z)),v\rangle=0$ for all $j=1$ to $d$, we have
$\langle\overline s_z\{C\}\phi(z,f(z)),v\wedge v^1\wedge
v^2\wedge...\wedge \widehat v^k\wedge...\wedge v^d\rangle=0$, as every
term of the expansion contains a factor of the form $\langle
\partial_wq_j(z,f(z)),v\rangle$.
Write the other term in (40) as $C(z)\sum_{i=1}^d R_i$ and consider
the $i^{th}$ term $R_i$ of this summation.  If $i\neq k$, $\langle
R_i,v\wedge v^1\wedge v^2\wedge...\wedge \widehat v^k\wedge...\wedge
v^d\rangle$ is zero because every term in the expansion contains a
factor of the form $\langle\sum_{j=1}^m{\partial q_k\over\partial
w_j}(z,f(z))dw_j,X\rangle$, where $X$ is one of
$v,v^1,v^2,v^3,...,v^{k-1},v^{k+1},...,v^d$, and every such factor is
zero by (41) and the definition of $v$.  Thus $\langle R,v\wedge
v^1\wedge v^2\wedge...\wedge\widehat v^k\wedge ...\wedge
v^d\rangle=C(z)\langle\sum_{i=1}^d R_i,v\wedge v^1\wedge
v^2\wedge...\wedge\widehat v^k\wedge ...\wedge v^d\rangle=C(z)\langle
R_k,v\wedge v^1\wedge v^2\wedge...\wedge\widehat v^k\wedge ...\wedge
v^d\rangle$.  Since (41) holds and $\langle\partial_w
q_j(z,f(z)),v\rangle=0$ for $j=1$ to $d$, the only nonzero term in the
expansion of $C(z)\langle R_k,v\wedge v^1\wedge v^2\wedge...\wedge
\widehat v^k\wedge...\wedge v^d\rangle$ is
$$\eqalign{&(-1)^{k-1}C(z)\langle (\sum_{j=1}^m{\partial
q_1\over\partial w_j}(z,f(z))dw_j),v^1\rangle\langle
(\sum_{j=1}^m{\partial q_2\over\partial
w_j}(z,f(z))dw_j),v^2\rangle...\cr
&\langle(\sum_{j=1}^m\sum_{\sigma=1}^\ell{\partial^2q_k\over\partial\overline
z_\sigma\partial w_j}(z,f(z))\overline s_z\{\overline z_\sigma\}dw_j+
\sum_{j,\sigma=1}^m{\partial^2q_k\over
\partial\overline w_\sigma\partial w_j}(z,f(z))\overline s_z \{\overline 
f_\sigma\} dw_j),v\rangle...\cr& \langle(\sum_{j=1}^m{\partial
q_d\over\partial w_j}(z,f(z))dw_j),v^d\rangle\cr}$$ which, by (41),
equals
$$\eqalign{&
(-1)^{k-1}C(z)\langle\sum_{j=1}^m\sum_{\sigma=1}^\ell{\partial^2q_k\over
\partial\overline
z_\sigma\partial w_j}(z,f(z))\overline s_z\{\overline
z_\sigma\}dw_j+\sum_{j,\sigma=1}^m{\partial^2q_k\over
\partial\overline w_\sigma\partial w_j}(z,f(z))\overline s_z \{\overline
f_\sigma\} dw_j,v\rangle\cr =&
(-1)^{k-1}C(z)\langle\overline\partial\partial q_k(z,f(z)),(\overline
{\sum_{j=1}^\ell s_z\{z_j\}{\partial\over\partial z_j}
+\sum_{j=1}^ms_z\{f_j\}{\partial\over\partial w_j}})\wedge
v\rangle.}$$ By (40) this quantity is zero, and this holds for $k=1$
to $d$.  Since $C(z)$ is never zero, we find that for all $k=1$ to
$d$,
$$\langle\overline\partial\partial
q_k(z,f(z)),(\overline {\sum_{j=1}^\ell
s_z\{z_j\}{\partial\over\partial z_j}
+\sum_{j=1}^ms_z\{f_j\}{\partial\over\partial w_j}})\wedge
v\rangle=0.$$  It is also true that for $i=1$ to $c$,
$$\langle\overline\partial\partial p_i(z,f(z)),(\overline
{\sum_{j=1}^\ell s_z\{z_j\}{\partial\over\partial z_j}
+\sum_{j=1}^ms_z\{f_j\}{\partial\over\partial w_j}})\wedge
v\rangle=0$$ since $p_i$ depends only on $z$ and $v$ has no terms
involving any ${\partial\over\partial z_j}$.  By definition of
$V(z,f(z))$ and since $v$ was chosen arbitrarily in $V(z,f(z))$, this
implies that
$$\sum_{i=1}^\ell s_z\{z_i\}{\partial\over\partial z_i}
+\sum_{i=1}^ms_z\{f_i\}{\partial\over\partial w_i}\in N(z,f(z)).$$
This is true for all $s_z\in \Nu^S_z$, so $\pi_{(z,f(z))}$
maps $N(z,f(z))$ onto $\Nu^S_z$.  We already know from Proposition 1
that this mapping is injective, so $\pi_{(z,f(z))}:N(z,f(z))\rightarrow
\Nu^S_z$ is an isomorphism, and $N(z,f(z))$ must have
the same complex dimension as $\Nu^S_z$, which is $\ell-c$.  The
foregoing argument holds for an arbitrary $z\in U$.  Since any point
$(z,w)$ in $M$ is on the graph of an $f$ defined on some such
$U\subset S$, we find that the complex dimension of $N(z,w)$ is
$\ell-c$ for all $(z,w)\in M$, as desired.

Now we have to prove (18).  Once again we fix $(z^0,w^0)\in M$ through
which there exists the graph of a CR map $f:U\rightarrow\Bbb C^m$ with
the properties in (III) for some neighborhood $U$ of $z^0$ in $S$.  We
note that every element of $N(z,f(z))$ for $z\in U$ has the form
$s_z+\sum_{i=1}^m s_z\{f_i\}{\partial\over \partial w_i}$ for some
$s_z\in\Nu^S_z$: if $n_z\in N(z,f(z))$ then
$\pi_{(z,f(z))}(n_z)\in\Nu^S_z$.  Let
$s_z=\pi_{(z,f(z))}(n_z)$.  Then $s_z+\sum_{i=1}^m
s_z\{f_i\}{\partial\over \partial w_i}\in N(z,f(z))$ and
$\pi_{(z,f(z))}( n_z)=\pi_{(z,f(z))}(s_z+\sum_{i=1}^m s_z\{f_i\}{\partial\over
\partial w_i})=s_z$.  By Proposition 1, $\pi_{(z,f(z))}:N(z,f(z))
\rightarrow \Nu^S_z$ is injective, so $n_z=s_z+\sum_{i=1}^m
s_z\{f_i\}{\partial\over \partial w_i}$, as claimed.

It suffices to prove (18) at points on the graph of a particular such
$f$, since the graphs of such $f$ foliate $M$ locally.  The dimension
of $N(z,w)$ is $\ell-c$ for all $(z,w)$ in $M$.  For all $z\in U$,
$N(z,f(z))\subset\Bbb CT_{(z,f(z))}M^f$, so the restriction of $N\oplus
\overline N$ to
$M^f$ is a subbundle of $\Bbb CTM^f$.  For a shrunken $U$, we may
construct the local basis for $\Bbb CTS$ near $z^0$ as in (19); we
write $G=U$.  Then we may construct the vector fields $T_i$ near
$(z^0,w^0)$ (see (20)) as at the beginning of the proof.  Since the
$T_i$ are commutators of vector fields whose values on $M^f$ are
tangent to $M^f$, the values of the $T_i$ are also tangent to $M^f$.
In fact, any commutator of vector fields in $\Gamma(G^M,N)$ and
$\Gamma(G^M,\overline N)$ has values on $M^f$ which are tangent to
$M^f$.


We conclude that to prove (18) for the vector fields $F_i$ indicated
there, it suffices to prove that if $z\in G$, $$\langle\partial
p_1(z,f(z))\wedge\partial p_2(z,f(z))\wedge...
\wedge\partial p_c(z,f(z))\wedge\partial q_k(z,f(z)),F^1_z
\wedge F^{2}_z\wedge...\wedge
F^{c+1}_z\rangle=0\eqno(42)$$ for all $F^i_z\in \Bbb
CT_{(z,f(z))}M^f$, $i=1$ to $c+1$. Now fix $z\in G$, so $(z,f(z))$ is
in the graph of $f$.  There is a complex linear isomorphism
$$I:\Bbb CT_0(\Bbb R^c\times\Bbb C^{\ell-c})\rightarrow \Bbb
CT_{(z,f(z))}M^f$$ which maps $\Bbb CT_0^{1,0}(\{0\}\times\Bbb
C^{\ell-c})$ onto $\Bbb CT^{1,0}_{(z,f(z))}M^f=N(z,f(z))$ and $\Bbb
CT_0^{0,1}(\{0\}\times\Bbb C^{\ell-c})$ onto $\Bbb
CT^{0,1}_{(z,f(z))}M^f=\overline N(z,f(z))$.  (Suppose $\Bbb R^c$ has
coordinates $x_i$, $i=1$ to $c$ and $\Bbb C^{\ell-c}$ has coordinates
$\zeta_i$, $i=1$ to $\ell-c$.  Then just let $I({\partial
\over\partial \zeta_i})=n_i(z,f(z))$,  $I({\partial
\over\partial \overline \zeta_i})=\overline n_i(z,f(z))$ and $I({\partial\over
\partial x_i})=T_i(z,f(z))$, where $n_i$ and $T_i$ were defined in 
(13) and (20), respectively.)  The map $I$ induces an isomorphism on
the cotangent spaces
$$I^*:\Bbb CT^*_{(z,f(z))}M^f\rightarrow\Bbb CT^*_0(\Bbb R^c\times\Bbb
C^{\ell-c}).$$ Note that all $\partial p_i(z,f(z))$ and $\partial
q_i(z,f(z))$ may be regarded as elements of $\Bbb CT^*_{(z,f(z))}M^f$,
by restriction of those forms to $\Bbb CT_{(z,f(z))}M^f$.  Let
$\psi_i$ be the cotangent in $\Bbb CT_0^*(\Bbb R^c\times\Bbb
C^{\ell-c})$ given by $\psi_i=I^*(\partial p_i(z,f(z)))$, for all
$i=1$ to $c$; more precisely, we have that if $\tilde F_z\in\Bbb
CT_0(\Bbb R^c\times\Bbb C^{\ell-c})$ then $\langle\psi_i,\tilde
F_z\rangle=\langle I^*(\partial p_i(z,f(z))),\tilde F_z\rangle
=\langle\partial p_i(z,f(z)),I(\tilde F_z)\rangle$.  Also write
$\xi_i=I^*(\partial q_i(z,f(z)))$ for all $i=1$ to $d$, so
$\langle\xi_i,\tilde F_z\rangle=\langle\partial q_i(z,f(z)),I(\tilde
F_z)\rangle.$ We know that if $\tilde F_z\in\Bbb
CT^{1,0}_0(\{0\}\times \Bbb C^{\ell-c})$ then $\langle\psi_i,\tilde
F_z\rangle=\langle\partial p_i(z,f(z)),I(\tilde F_z) \rangle=0$ for
all $i=1$ to $c$ since $I(\tilde F_z)$ is a $(1,0)$ tangent to $M^f$
at $(z,f(z))$.  Similarly, $\langle\xi_i,\tilde F_z\rangle=0$ for all
$i=1$ to $d$.  We also have that $\langle\psi_i,\tilde F_z\rangle=0$
and $\langle\xi_i,\tilde F_z\rangle=0$ for all $\psi_i,\xi_i$ and
$\tilde F_z\in \Bbb CT^{0,1}_0(\{0\}\times \Bbb C^{\ell-c})$ since
then $I(\tilde F_z)$ is a $(0,1)$ tangent to $M^f$ at $(z,f(z))$.  By
Lemma 2,
$\psi_1\wedge\psi_2\wedge\psi_3\wedge...\wedge\psi_c\wedge\xi_k=0$ as
a $c+1$-cotangent in $\Bbb CT_0^*(\Bbb R^c\times\Bbb C^{\ell-c})$, so
for all $\tilde F^1_z,\tilde F^2_z,...,\tilde F^{c+1}_z\in\Bbb
CT_0(\Bbb R^c\times\Bbb C^{\ell-c})$,
$$\langle\psi_1\wedge\psi_2\wedge\psi_3\wedge...\wedge\psi_c\wedge\xi_k,
\tilde F^1_z,\wedge\tilde F^2_z\wedge...\wedge \tilde F^{c+1}_z\rangle=0.
\eqno(43)$$
Temporarily writing $\psi_{c+1}=\xi_k$, we find from (43) that
$$0=\sum_\sigma\hbox{\rm sgn}(\sigma)\prod_{i=1}^{c+1}\langle\psi_i,
\tilde F^{\sigma(i)}_z\rangle =\sum_\sigma\hbox{\rm
sgn}(\sigma)\langle\partial q_k(z,f(z)),I(\tilde
F^{\sigma(c+1)}_z)\rangle\prod_{i=1}^{c}\langle\partial
p_i(z,f(z)),I(\tilde F^{\sigma(i)}_z)\rangle,\eqno(44)$$ where
$\sigma$ ranges over all permutations of $\{1,2,3,...,c+1\}$ and
$\hbox{\rm sgn}(\sigma)$ is the sign of such a permutation $\sigma$.
Then (44) implies that
$$\langle\partial
p_1(z,f(z))\wedge\partial p_2(z,f(z))\wedge...
\wedge\partial p_c(z,f(z))\wedge\partial q_k(z,f(z)),I(\tilde
F^1_z)\wedge I(\tilde F^2_z)\wedge...\wedge I(\tilde
F^{c+1}_z)\rangle=0\eqno(45)$$ for all $\tilde F^1_z,\tilde
F^2_z,...,\tilde F^{c+1}_z\in\Bbb CT_0(\Bbb R^c\times\Bbb
C^{\ell-c})$.  Since $I$ is an isomorphism, as the $\tilde F^i_z$
range over $\Bbb CT_0(\Bbb R^c\times\Bbb C^{\ell-c})$, $I(\tilde
F^i_z)$ ranges over $\Bbb CT_{(z,f(z))}M^f$.  Thus we conclude from
(45) that (42) holds for all $F^i_z\in \Bbb CT_{(z,f(z))}M^f$, $i=1$
to $c+1$.  As observed earlier, this implies (18).  That completes the
proof that (III) implies (IV).

Now we must prove that if (I),(II),(III),(IV) hold then the last
statements of the theorem hold.  While proving that (II) implies
(III), we proved that the graphs in (III) arise as integral manifolds
of $\Re\,{\Cal T}$ and that the $(1,0)$ tangent space to such a graph
$f$ at $(z,f(z))$ is $N(z,f(z))$.  Lastly, suppose a CR map $g$ exists
as stated in the theorem.  Let $M^g$ denote its graph.  Then the
complexified tangent space $\Bbb CT_{(z,g(z))}M^g$ to the graph of $g$
contains $N(z,g(z))$ (as assumed) so it contains
$N(z,g(z))\oplus\overline N(z,g(z))$ (since the graph of $g$ is a real
manifold, $\Bbb CTM^g$ is self-conjugate) so it contains ${\Cal
T}(z,g(z))$ (since $\Bbb CTM^g$ is involutive, $\Bbb CTM^g_{(z,g(z))}$
must contain $T_i(z,g(z))$ for $i=1$ to $c$, where $T_i$ is defined in
(20)).  In fact, for all $z$ in the domain of $g$, $\Bbb
CT_{(z,g(z))}M^g$ equals ${\Cal T}(z,g(z))$ because each has
complex dimension $2\ell-c$.  Thus the graph of $g$ is an integral
manifold for $\Re\,{\Cal T}$, so by the Frobenius theorem it must be
one of the graphs from (III).$\hfill\square$

The next theorem is a statement about the existence of CR maps
whose graphs lie in $M$ and whose domains are all of $S$.

\proclaim{Theorem 2}  Suppose that $S,M$ satisfy (1),(2), respectively,
and that $S$ is connected, simply connected and of finite
type $\tau$ at every point.  Suppose also that for every compact $K
\subset S$, $M_K\equiv\{(z,w)\in M:z\in K\}$ is compact.
Suppose that any of properties (I),(II),(III) or (IV) of Theorem 1
hold.  Then for every $(z^0,w^0)\in M$ there exists a unique CR map
$f:S\rightarrow\Bbb C^m$ whose graph is an integral manifold of
$\Re\,{\Cal T}$ and which passes through $(z^0,w^0)$.  Let
$\phi_{i_1,i_2,...,i_d}(z,w)$ be the function defined in (III) of
Theorem 1.  If there exist CR functions $h_{i_1,i_2,...,i_d}(z)$
defined for $z\in S$ such that
$$C(z)\equiv\sum_{1\leq i_1<i_2<...<i_d\leq m} h_{i_1,i_2,...,i_d}(z)
\phi_{i_1,i_2,...,i_d}(z,f(z))\neq 0\eqno(46)$$ for all $z\in S$
then $$\hbox{$z\mapsto{1\over C(z)}\phi_{i_1,i_2,...,i_d}(z,f(z))$ is
a CR function on $S$ for integers $i_j$,}\eqno(47)$$
$1\leq i_1<i_2<...<i_d\leq m$.
If $d=m$ then such $h_{i_1,i_2,...,i_d}$ will exist.
\endproclaim 

\noindent{\bf Note:}
Theorems 3 and 4 give natural circumstances where functions
$h_{i_1,i_2,...,i_d}$ exist such that property (46) holds.

\noindent{\bf Proof:}
By Theorem 1, properties (I),(II),(III) and (IV) are all equivalent.
Thus there exists an involutive subbundle $\Re\,{\Cal T}$ of the real
tangent bundle to $M$ whose integral manifolds are locally graphs of
CR maps on open subsets of $S$.  Pick $(z^0,w^0)\in M$ and let $L$ be
the maximal connected integral manifold of $\Re\,{\Cal T}$ which
passes through $(z^0,w^0)$.  (For existence of such a manifold, see
[Wa, Theorem 1.64].)  We claim that $L$ is the graph of a CR map
$f:S\rightarrow\Bbb C^m$.  ($L$ is the union of submanifolds which are
graphs over open subsets of $S$; the topology we use for $L$ is that
generated by the topologies of these graphs taken together.)  First we
claim that the projection mapping $P$ from $L$ to $S$ is a covering
map.  (We use the definition from [Ma, p. 118].)  We show that it is
surjective and to this end we show that $P(L)$ is open in $S$.  Fix
$z\in S$ which is in $P(L)$; then there exists $(z,w)\in L$ and an
integral manifold of $\Re\,{\Cal T}$ passing through $(z,w)$ which is
the graph of a CR function in a neighborhood $U$ of $z$.  This graph
is contained in $L$ since the maximal integral manifold $L$ of
$\Re\,{\Cal T}$ contains any integral manifold of $\Re\,{\Cal T}$
passing through a point of $L$.  (See [Wa, Theorem 1.64].) Thus $P(L)$
contains $U$.  This shows that the projection of $L$ to $S$ is open in
$S$.  If $P(L)$ is not all of $S$, then since $S$ is connected and
$P(L)$ is open, there exists a $C^0$ path $\gamma:[0,1]
\rightarrow S$ such that $\gamma([0,1))\subset P(L)$, $\gamma(0)=z^0$,
but $\gamma(1)\notin P(L)$.  We claim that there exists a continuous
$\gamma': [0,1]\rightarrow L$ such that $\gamma'(0)=(z^0,w^0)$ and
$P\circ\gamma'=\gamma$.

Before proving this, suppose we have two continuous paths
$\gamma',\gamma''$ such that $\gamma':[0,\epsilon']
\rightarrow L$, $\gamma'':[0,\epsilon'']
\rightarrow L$, $\gamma'(0)=\gamma''(0)=(z^0,w^0)$, $\epsilon'\leq\epsilon''$,
and $P\circ\gamma'=P\circ\gamma''=\gamma$ on $[0,\epsilon']$.  Then in
fact $\gamma'=\gamma''$ on $[0,\epsilon']$.  To see this, note that
$\gamma'(0)=\gamma''(0)$, so the domain of coincidence of
$\gamma',\gamma''$ is nonempty.  The set where $\gamma'=\gamma''$ is
closed in $[0,\epsilon']$ because $\gamma',\gamma''$ are continuous.
It is also open:if $t_1$ is a point where $\gamma'(t_1)=\gamma''(t_1)$
then in a neighborhood of that point in $L$, $L$ is a graph over an
open neighborhood of $P(\gamma'(t_1))$ in $S$, so $P$ maps a
neighborhood of $\gamma'(t_1)$ in $L$ homeomorphically to a
neighborhood of $P(\gamma'(t_1))$ in $S$.  For $t$ near $t_1$,
$\gamma'(t),\gamma''(t)$ lie in that neighborhood in $L$ of
$\gamma'(t_1)$, so $\gamma'(t)=\gamma''(t)$ for $t$ near $t_1$.  Thus
the set where $\gamma'=\gamma''$ in $[0,\epsilon']$ is open, closed
and nonempty, so equals $[0,\epsilon']$, as desired.

Now let $\gamma'(0)=(z^0,w^0)$.  Near $\gamma'(0)$ in $L$, $L$ is a
graph over a neighborhood of $\gamma(0)$; thus for some $\epsilon>0$
we may define a continuous $\gamma': [0,\epsilon]\rightarrow L$ such
that $P\circ\gamma'=\gamma$ on $[0,\epsilon]$.  Let $t_0$ be the
supremum of all $t\in [0,1]$ such that we may define a continuous
$\gamma' : [0,t]\rightarrow L$ such that $\gamma'(0)=(z^0,w^0)$ and
$P\circ\gamma' =\gamma$.  Any two such paths coincide on their common
domains by the previous paragraph.  Then there exists a continuous
$\gamma' : [0,t_0)\rightarrow L$ such that $P\circ\gamma' =\gamma$.
Let $K=\gamma([0,t_0])$ and let $(z,w)$ be a limit point in $M_K$ of
the open path $\gamma' ([0,t_0))$ of the form
$\lim_{n\rightarrow\infty}\gamma'(t_n)$ where $t_n\uparrow t_0$.
(Recall that $M_K$ is compact.)  By Theorem 1, there exists a
neighborhood $U^M$ of $(z,w)$ in $M$ which is foliated by graphs of CR
functions on an open subset $U$ of $S$ where the graphs are integral
manifolds of $\Re\,{\Cal T}$.  For $t$ near $t_0$, $\gamma(t)$ lies in
$U$ and for large $n$, $\gamma' (t_n)$ lies in $U^M$ on one of the
foliating graphs $f^n:U\rightarrow\Bbb C^m$.  Since the image of
$\gamma' $ is connected and $\gamma(t)$ lies in $U$ for $t$ near
$t_0$, $\gamma'(t)$ lies only on the graph of a particular CR
$f:U\rightarrow\Bbb C^m$ for $t$ near $t_0$.  The graph of $f$ is an
integral manifold of $\Re\,{\Cal T}$ for some points in $L$ (the
points on the path $\gamma' (t)$ for $t$ near $t_0$.)  Since $L$ is
maximal, $L$ contains the graph of $f$.  Furthermore, we may use this
fact to extend $\gamma' $ to a neighborhood of $t_0$ in $[0,1]$.  If
$t_0<1$, this contradicts the maximality of $t_0$, so $t_0=1$ and we
have the path $\gamma' $ as desired.  However, this implies that
$P(L)$ contains $\gamma(1)$, a contradiction, so the assumption that
$P(L)$ is not all of $S$ is false: we have $P(L)=S$.


We recall again from Theorem 1 that given any $(z,w)\in M$ there
exists an open neighborhood $U$ of $z$ in $S$ and an open neighborhood
of $U^M$ of $(z,w)$ in $M$ such that $U^M$ is the disjoint union of
graphs of CR maps $f:U\rightarrow\Bbb C^m$, where the graphs are all
integral manifolds of $\Re\,{\Cal T}$.  Fixing $z=z'$, we find such
open sets $U_w,U^M_w$ for every $w\in\Bbb C^m$ such that $(z',w)\in
M$.  Since $M_{z'}$ is compact, finitely many of the $U^M_w$ cover
$M\cap\{(z,w)\in M:z=z'\}$.  Take the (finite) intersection of all
associated $U_w$ and let $U$ be the path component of it which
contains $z'$.  Then $U$ is open (as a path component of the finite
intersection of open sets) and given any point of the form $(z',w)\in
M$ there exists a CR map $f:U\rightarrow\Bbb C^m$ whose graph contains
$(z',w)$ and which is an integral manifold for $\Re\,{\Cal T}$.
(I.e., given fixed $z'\in S$ there exists a fixed $U$ for all $w\in
M_{z'}$ such that such an $f$ exists.)  Thus given $z'\in S$ we have
found a path connected neighborhood $U$ of $z'$ in $S$ such that
$P^{-1}(U)$ is the union of disjoint CR graphs over $U$ each of which
is a path connected open subset of $L$ and each of which projects onto
$U$ through $P$.  This proves that $P:L\rightarrow S$ is a covering
map.

Next we claim that $P:L\rightarrow S$ is injective; this will show
that $L$ is a graph.  Suppose it is not injective; then there exist
$z^0\in S$, $(z^0,w^0)\in L$ and $(z^0,w^1)\in L$ with $w^0\neq w^1$.
Since $L$ is a connected manifold, it is path connected, so there
exists a path in $L$ from $(z^0,w^0)$ to $(z^0,w^1)$ which projects to
a path in $S$ from $z^0$ to $z^0$.  By the path lifting lemma (see
[Ma, Lemma 3.3]), since $S$ is simply connected and $L$ a covering
space of $S$, we have $(z^0,w^0)=(z^0,w^1)$, as desired.  Thus
$P:L\rightarrow S$ is injective, as desired and $L$ is the graph of a
mapping $f:S\rightarrow\Bbb C^m$.  Then $f$ is CR because $L$ is the
union of graphs of CR mappings defined on open subsets of $S$.  If
another mapping $\tilde f:S\rightarrow\Bbb C^m$ exists with the
properties $f$ has in Theorem 2, then its graph is an integral
manifold of $\Re\,{\Cal T}$ which passes through $(z^0,w^0)$.  Since
the graph of $f$ is $L$, the maximal connected integral manifold of
$\Re\,{\Cal T}$ which passes through $(z^0,w^0)$, the graph of $\tilde
f$ is contained in the graph of $f$, so $f=\tilde f$, as desired.

We need to show that $z\mapsto{1\over
C(z)}\phi_{j_1,j_2,...,j_d}(z,f(z))$ is CR; if $\overline
s_z\in\overline \Nu^S_z$, then for all integers $j_1,j_2,...,j_m$ such
that $1\leq j_1<j_2<...<j_d\leq m$ we have
$$\eqalign{&\overline s_z{\phi_{j_1,j_2,...,j_d}(\cdot,f(\cdot))\over
C(\cdot)}\cr=&{1\over C(z)^2}
\biggl(\sum_{1\leq
i_1<i_2<...<i_d\leq m} h_{i_1,i_2,...,i_d}(z)
\phi_{i_1,i_2,...,i_d}(z,f(z))\overline s_z\{\phi_{j_1,j_2,...,j_d}
(\cdot,f(\cdot))\}\cr&-\phi_{j_1,j_2,...,j_d}(z,f(z))
\sum_{1\leq
i_1<i_2<...<i_d\leq m} h_{i_1,i_2,...,i_d}(z)
\overline s_z\{\phi_{i_1,i_2,...,i_d}(\cdot,f(\cdot))\}\biggr)\cr
=&{1\over C(z)^2}\sum_{1\leq i_1<i_2<...<i_d\leq m}
h_{i_1,i_2,...,i_d}(z)\biggl(\overline
s_z\{\phi_{j_1,j_2,...,j_d}(\cdot,f(\cdot))\}\phi_{i_1,i_2,...,i_d}(z,f(z))\cr
&-\overline
s_z\{\phi_{i_1,i_2,...,i_d}(\cdot,f(\cdot))\}
\phi_{j_1,j_2,...,j_d}(z,f(z))\biggr)
\cr
=&0}$$ from Lemma 4.  (Note that $\overline
s_z\{h_{i_1,i_2,...,i_d}\}=0$ for all $z\in S$ since
$h_{i_1,i_2,...,i_d}$ is CR.)  This being true for all $\overline
s_z\in\overline
\Nu^S_z$ and all $z\in S$, we find that $z\mapsto{1\over
C(z)}\phi_{j_1,j_2,...,j_d}(z,f(z))$ is CR as desired.

If $d=m$ then the note before the proof of Theorem 1 applies: there is
only one $\phi_{i_1,i_2,...,i_d}=\phi_{1,2,3,...,m}$ and we may just let
$C(z)=1/\phi_{1,2,3,...,m}(z,f(z))$ which is never zero for $z\in S$.
$\hfill\square$

If the conditions of Theorem 2 hold then we let ${\Cal F}$ be the set
of CR mappings $f$ whose graphs lie in $M$ and are maximal integral
manifolds of $\Re{\Cal T}$.

\proclaim{Corollary 1}  Suppose that $S$ and $M$ satisfy the requirements
of Theorem 2 and that in addition $S$ is the boundary of a bounded
domain $D$ in $\Bbb C^\ell$, $\ell\geq 2$.  Then the CR maps in ${\Cal F}$
all extend to be continuous on $\overline D$ and analytic on $D$;
these extensions are solutions to (RH).

\endproclaim

\noindent{\bf Proof:}
The CR maps which arise from Theorem 2 are $C^\infty$ on $S$; they
extend to be analytic on $D$ by the global CR extension theorem.
(This is the theorem commonly known as Bochner's extension theorem; we
choose not to use this attribution becuase of the conclusions in the
paper [Ra].)  \hfill $\square$

Corollary $1$ provides circumstances where the Riemann-Hilbert problem
(RH) for $M$ is solvable.  If $f\in{\Cal F}$ then we let $\hat f=
(\hat f_1,\hat f_2,...,\hat f_m)$ denote the analytic extension of $f$
to $D$, if the extension exists.

\proclaim{Theorem 3}   Suppose that $S,M$ satisfy (1),(2), respectively,
and that $S$ is connected, simply connected and of finite type $\tau$
at every point.  Suppose also that $M_K$ is compact for every compact
$K$ in $S$ and that for all $z\in S$, $M_z$ is a convex hypersurface
which encloses the origin in $\Bbb C^m$.  Lastly suppose that any of
properties (I),(II), (III) or (IV) of Theorem 1 hold.  Then $M$ is the
disjoint union of CR maps $f:S\rightarrow
\Bbb C^m$ such that there exists a nonzero complex-valued function 
$C(z)$ defined for $z\in S$ for which
$$z\mapsto{1\over C(z)}{\partial q_1\over\partial
w_i}(z,f(z))\eqno(48)$$ is a CR function on $S$ for $i=1,2,3,...,m$.
\endproclaim

\noindent{\bf Proof:}
From Theorem 2, we may conclude that for any point in $M$ there exists
a CR mapping $f:S\rightarrow\Bbb C^m$ whose graph is contained in $M$
and passes through that point.  Since (48) is nothing more than
condition (47) in the case $d=1$, we must only show the existence of
$h_i$, $i=1$ to $m$, such that (46) holds in the case $d=1$, i.e. that
$$\sum_{i=1}^mh_i(z){\partial q_1\over\partial w_i}(z,f(z))\neq
0\eqno(49)$$ for all $z\in S$, where $h=(h_1,h_2,...,h_m)$.  If
$f=(f_1,f_2,...,f_m)$ then it will suffice to let $h_i=f_i$ for $i=1$
to $m$.  The reason for this is that by convexity of the surfaces
$M_z$ for all $z\in S$, the complex tangent plane to $M_z$ in $\Bbb
C^m$ at $f(z)$ does not pass into the region enclosed by $M_z$, so
does not pass through the origin of $\Bbb C^m$.  This complex tangent
plane has the form
$$\{v=(v_1,v_2,v_3,...,v_m)\in\Bbb C^m:\sum_{i=1}^m {\partial
q_1\over\partial w_i}(z,f(z)) v_i=\zeta_z\}$$ for some complex constant
$\zeta_z$.  We cannot have $\zeta_z=0$ because the plane does not pass
through the origin.  Thus $\sum_{i=1}^m {\partial q_1\over\partial
w_i}(z,f(z)) f_i(z)=\zeta_z\neq 0$ for all $z\in S$ since $f(z)$
belongs to the tangent plane to $M_z$ at $f(z)$.  Thus (49) holds and
Theorem 3 holds.$\hfill\square$

{\bf Note:} Since all that is required in Corollary 2 is that
the {\it complex} tangent planes to $M_z$ not pass into the interior
of the region enclosed by $M_z$, the theorem need only require that $M_z$
enclose a region that is {\it lineally convex} or {\it hypoconvex}.
(See [Ki], [Wh] or [H\" o, p. 290].)

\proclaim {Corollary 2}  Suppose $S$ and $M$ satisfy all the conditions
of Theorem 3 and that in addition $S$ is the boundary of a bounded
domain $D$ in $\Bbb C^\ell$, $\ell\geq 2$.  Then the CR maps in ${\Cal
F}$ which arise from Theorem 3 all extend to be continuous on
$\overline D$ and analytic on $D$; these extensions are solutions to
(RH).
\endproclaim

\noindent{\bf Proof:}
This is again an application of the global CR extension
theorem.$\hfill\square$
 
If $Y$ is any compact subset of $\Bbb C^n$ then the {\it polynomial
hull} of $Y$ is the set
$$\hbox {$\widehat Y =\{z\in\Bbb C^n\bigl ||P(z)|\leq {\displaystyle
\sup _{w\in Y}} |P(w)|$ for all polynomials $P$ on $\Bbb C^n\}$.}$$ We
say that $Y$ is {\it polynomially convex} if $\widehat Y=Y$.  Theorem
4 shows that under some circumstances, the properties that $f$
possesses in (III) of Theorem 1 guarantee that $f$ possesses an
extremal property.

\proclaim{Theorem 4}  Suppose that $S$,$M$ satisfy (1),(2) and
that $S$ is a hypersurface which bounds a bounded strictly
pseudoconvex open set $D$ such that $\overline D$ is polynomially
convex.  Further suppose that $S$ is connected and simply connected.
Suppose that in (2), defining function $q_d(z,w)$ satisfies the
property that for all $z\in S$, $q_d(z,w)$ is strictly convex as a
function of $w$ in a neighborhood of $M_z$.  Suppose that for
$i=1,2,3,...,d-1$, defining function $q_i$ has the form
$$q_i(z,w)=\Re\left(\sum_{j=1}^m\alpha ^i_j(z)w_j\right),$$ for some
matrix $(\alpha ^i_j)_{i,j=1}^m$ of functions analytic in a
neighborhood of $\overline D$, where the determinant of $(\alpha
^i_j)$ is never zero on $\overline D$.  Let $\tilde M=\{(z,w)\in
\overline D\times\Bbb C^m:q_i(z,w)=0,i=1,2,...,d-1\}$ and $\tilde M
_z=\{w:(z,w)\in\tilde M\}$.  Assume that $M$ is compact and that for
all $z\in S$ the origin of $\Bbb C^m$ is in the bounded (convex)
component of $\tilde M_z\setminus M_z$.  Suppose that any of the
properties (I),(II), (III) or (IV) of Theorem 1 hold.  Then the set
${\Cal F}$ is well defined and for all $f\in{\Cal F}$, the graph of
$\hat f$ in $\overline D\times\Bbb C^m$ lies in the boundary of
$\widehat M$ as a subset of $\tilde M$.  In particular, given
$f\in{\Cal F}$ and some $z^0\in D$, the only continuous mapping
$k:\overline D\rightarrow\Bbb C^m$ such that $k$ is analytic in $D$,
$k(z^0)=\hat f(z^0)$ and $k(z)$ belongs to the convex hull of $M_z$
for all $z\in S$ is $k=\hat f$.
\endproclaim
{\bf Note:} By Proposition 1(iii), in order to verify that (I) of Theorem 1
holds we must only check that (17) holds.

\noindent{\bf Proof:}  Since $D$ is strictly pseudoconvex,
$S$ is of type $2$.  Thus the conditions of Theorem 2 and Corollary 1
are satisfied, so ${\Cal F}$ exists.  By an analytic change of
variable which is linear in $w$, we may assume without loss of
generality that $q_i(z,w)=2\,\Re\, w_i$ for $i=1,2,3,...,d-1$.  Fix
$f\in{\Cal F}$.  Then $\Re\,f_i(z)=0$ for $i=1,2,...,d-1$ and for
$z\in S=\partial D$, so $\Re\,f_i(z)=0$ for $z\in\overline D$ and the
graph of $\hat f$ over $\overline D$ is contained in $\tilde M$.
Then, letting $\phi$ denote the $d$-form defined in Theorem 1, $$\phi
(z,w)=\sum_{j=d}^m {\partial q_d\over\partial w_j}(z,w)dw_1\wedge
dw_2\wedge dw_3\wedge...\wedge dw_{d-1}\wedge dw_j,$$ where
$\phi_{1,2,3,...,d-1,j}(z,w)={\partial q_d\over\partial w_j}(z,w)$,
$j=d$ to $m$.  Furthermore, $\sum_{j=d}^m{\partial q_d\over\partial
w_j}(z,f(z))f_j(z)$ is nonzero for all $z\in S$ because of the fact
that $M_z$ is a convex hypersurface in $\{w\in\Bbb C^m:\Re\,w_i=0,
i=1,2,3,...,d-1\}$ such that $M_z$ encloses the origin.  (Reasoning is
similar to that used in Theorem 3.)  This means that (46) holds, where
$h_{1,2,3,...,d-1,j} =f_j$ for $j=d$ to $m$, so $C(z)=\sum_{j=d}^m
f_j(z){\partial q_d\over\partial w_j}(z,f(z))$.  Thus conclusion (47)
of Theorem 2 holds.  Let $g_j(z)={1\over C(z)}{\partial
q_d\over\partial w_j}(z,f(z))$ for $z\in S$ and $j=d,d+1,...,m$.  For
$j=1$ to $d-1$ let $g_j$ be the zero function on $\overline D$.  Then
by (47), $g_j$ is CR on $S$.  Let $g:S\rightarrow
\Bbb C^m$ be given by $g=(g_1,g_2,...,g_m)=(0,0,...,0,g_d,
g_{d+1},...,g_m)$.  Then $g$ extends continuously to $\overline D$ and
analytically to $D$ to produce a mapping $\hat g=(\hat g_1,\hat
g_2,...,\hat g_m)$, by the global CR extension theorem.  Since
$\sum_{i=1}^mf_i(z)g_i(z)=1$ for all $z\in S$ and both $f$ and $g$
extend analytically to $D$, neither $\hat f$ nor $\hat g$ can be zero
anywhere on $\overline D$.  Let $B_m(R)$ be the open ball of radius
$R$ about the origin in $\Bbb C^m$.  Since $\overline D$ is
polynomially convex, so is $\overline D\times \overline {B_m(R)}$ for
any $R>0$.  Choose $R>0$ so large that $\widehat M\subset\overline
D\times \overline {B_m(R)}$.  Since $\hat g$ is never zero on
$\overline D$, the set
$$W(z,t)\equiv\{w\in\Bbb C^m:\sum_{i=1}^m \hat g_i(z)w_i=t\}$$ is a
complex hyperplane in $\Bbb C^m$ for any $z\in \overline D$ and $t>0$.
For any $z\in \overline D$ and $t>1$, this hyperplane comes no closer
than $t/K$ to the origin, where $K$ is the maximum modulus of $\hat g$ on
$\overline D$.  Fix $t_0$ so large that $t_0/K$ is greater than $R$.
Let 
$$W(t)\equiv\{(z,w)\in\overline D\times\Bbb C^m:\sum_{i=1}^m
\hat g_i(z)w_i=t\}.$$ Then if $t=t_0$, $W(t_0)$ does not meet $\widehat M$.
Let $t_1$ be the infimum of all $t$ such that $W(t)$ does not meet
$\widehat M$ and suppose that $t_1>1$.  Then the function
$F_t(z,w)=1/(-t+(\sum_{i=1}^m \hat g_i(z)w_i))$ is defined in a
neighborhood of $\widehat M$ in $\overline D\times\Bbb C^m$ for
$t_1<t\leq t_0$.  Since $D$ is strictly pseudoconvex, $\hat g_i$ is
uniformly approximable by functions analytic in a neighborhood of
$\overline D$.  (This theorem comes from [Li,He].)  Since $\overline
D$ is polynomially convex, the Oka-Weil Theorem in turn guarantees
that $\hat g_i$ is uniformly approximable by polynomials on $\overline
D$ for $i=1,2,...,m$, so $F_t$ is the uniform limit on $\widehat M$ of
functions analytic in a neighborhood of $\widehat M$ in $\Bbb
C^\ell\times\Bbb C^m$ for $t_1<t\leq t_0$.  By the Oka-Weil theorem,
for $t_1<t\leq t_0$, $F_t$ is uniformly approximable by polynomials on
$\widehat M$, so the supremum of $F_t$ on $\widehat M$ is less than or
equal to the supremum of $F_t$ on $M$, by the definition of polynomial
hull.  For the same $t$, we can see that $F_t$ is uniformly bounded on
$M$: for every $z\in S$, $M_z$ is a strictly convex hypersurface in
$\tilde M_z\equiv\{w\in\Bbb C^m:(z,w)\in\tilde M\}$, and $W(z,1)\cap
\tilde M_z$ is
tangent to $M_z$ in $\tilde M_z$ at $f(z)$.  Thus for all $t>1$,
$W(z,t)$ is a dilation of $W(z,1)$ away from the origin of $\Bbb C^m$
(which is contained in the convex hull of $M_z$), so $W(z,t)$ is
disjoint from $M_z$.  For $t>t_1>1$, the distance from $W(z,t)$ to
$M_z$ is then bounded below in $z\in
\partial D$, so $F_t$ is uniformly bounded on $M$ for $t>t_1>1$.
We conclude, from the observation that the $F_t$ are uniformly
approximable by polynomials on $\widehat M$ for $t_1<t<t_0$, that the
$F_t$ are uniformly bounded on $\widehat M$ for $t_1<t\leq t_0$.  This
is impossible because for $t$ near $t_1$, the singularity set $W(t)$
of $F_t$ approaches $\widehat M$ by definition of $t_1$.  Thus we must
have $t_1=1$ and $W(t)$ is external to $\widehat M$ for $t>1$.  Note
that the graph of $z\mapsto t\hat f(z)$ lies in $\tilde M$ and lies in
$W(t)$ since $\sum_{i=1}^m \hat g_i(z)\hat f_i(z)=1$ for all
$z\in\overline D$.  Thus the graph of $t\hat f$ lies external to
$\widehat M$ but inside $\tilde M$.  This being true for all $t>1$, we
find that the graph of $\hat f$ is in the boundary of $\widehat M$ as
a subset of $\tilde M$ since every point $(z,\hat f(z))$ on the graph
of $\hat f$ over $\overline D$ is the limit point of a set of points
$\{(z,t\hat f(z)):t>1\}$ in $\tilde M$ which is external to $\widehat
M$.

Lastly, suppose $k$ exists as indicated in the theorem.  Then the
function $z\mapsto \sum_{i=1}^m \hat g_i(z)k_i(z)-t$ is nonzero for
$t>1$ and $z\in\partial D$ since (as we showed) $W(z,t)$ is disjoint
from $M_z$ (so disjoint from the convex hull of $M_z$ as well) for all
$z\in\partial D$, $t>1$.  Thus $z\mapsto \sum_{i=1}^m \hat
g_i(z)k_i(z)-t$ is nonzero for $z\in D$ also.  If we take a limit as
$t\rightarrow 1^+$ we find from Hurwitz' theorem that $z\mapsto
\sum_{i=1}^m \hat g_i(z)k_i(z)-1$ is either identically zero on $D$ or
nonzero for $z\in D$.  The latter cannot be the case because
$\sum_{i=1}^m \hat g_i(z^0)k_i(z^0)-1=\sum_{i=1}^m \hat g_i(z^0)\hat
f_i(z^0)-1=0$.  Since $W(z,1)\cap\tilde M_z$ is a tangent plane to
$M_z$ in $\tilde M_z$ at $f(z)$ for all $z\in
\partial D$, the only point $w$ on both $W(z,1)$ and the convex hull of $M_z$ 
is $f(z)$.  But since $\sum_{i=1}^m \hat g_i(z)k_i(z)-1=0$ for all
$z\in\partial D$, $k(z)$ is in $W(z,1)$; it is in the convex hull of
$M_z$ by assumption.  Thus $k(z)=\hat f(z)$ for all $z\in\partial D$,
so for all $z\in \overline D$, as desired.  $\hfill\square$

In the case when $M$ has fibers $M_z$ which are real hypersurfaces
in $\Bbb C^m$ (i.e., $d=1$), the statement of Theorem 4 is of
interest:

\proclaim{Corollary 3}  Suppose that $S$,$M$ satisfy (1),(2) and
that $S$ is a hypersurface which bounds a bounded strictly
pseudoconvex open set $D$ such that $\overline D$ is polynomially
convex.  Further suppose that $S$ is connected and simply connected.
Suppose that $M$ is compact and that for every $z\in S$, $M_z$ is a
hypersurface which encloses a strictly convex open set in $\Bbb C^m$
such that the origin of $\Bbb C^m$ lies in that open set.  Suppose
that any of the properties (I),(II), (III) or (IV) of Theorem 1 hold.
Then the set ${\Cal F}$ is well defined and for all $f\in{\Cal F}$,
the graph of $\hat f$ in $\overline D\times\Bbb C^m$ lies in the
boundary of $\widehat M$.  In particular, given $f\in{\Cal F}$ and
some $z^0\in D$, the only mapping $k:\overline D\rightarrow\Bbb C^m$
such that $k$ is analytic in $D$, $k(z^0)=\hat f(z^0)$ and $k(z)$
belongs to the convex hull of $M_z$ for all $z\in S$ is $k=\hat f$.
\endproclaim

\noindent{\bf Proof:}  
This is simply the case $d=1$ of Theorem 4. $\hfill\square$

In Corollary 4, we consider the case where the functions $\alpha^i_j$
are only defined for $i=1$ to $d-1$.

\proclaim{Corollary 4} Suppose that $S$,$M$ satisfy (1),(2) and
that $S$ is a hypersurface which bounds a strictly pseudoconvex
bounded open set $D$ such that $\overline D$ is polynomially convex.
Further suppose that $S$ is connected and simply connected.  Suppose
that in (2), defining function $q_d(z,w)$ satisfies the property that
for all $z\in S$, $q_d(z,w)$ is strictly convex as a function of $w$
in a neighborhood of $M_z$.  Suppose that for $i=1,2,3,...,d-1$,
defining function $q_i$ has the form
$$q_i(z,w)=\Re\left(\sum_{j=1}^m\alpha^i_j(z)w_j\right),$$ for some
matrix $(\alpha^i_j){_{j=1}^m}{_{i=1}^{d-1}}$ of functions analytic in
a neighborhood of $\overline D$.  Let $\tilde M$ be as before.  Assume
that $M$ is compact and that the origin of $\Bbb C^m$ is in the
bounded (convex) component of $\tilde M_z\setminus M_z$ for
all $z\in S$.  Suppose that any of the properties (I),(II), (III) or
(IV) of Theorem 1 hold, that $m\geq 2\ell$ and that $d-1\leq
\ell$.  Then the set
${\Cal F}$ is well defined and for all $f\in{\Cal F}$, the graph of
$\hat f$ in $\overline D\times\Bbb C^m$ lies in the boundary of
$\widehat M$ as a subset of $\tilde M$.  In particular, given
$f\in{\Cal F}$ and some $z^0\in D$, the only continuous mapping
$k:\overline D\rightarrow\Bbb C^m$ such that $k$ is analytic in $D$,
$k(z^0)=\hat f(z^0)$ and $k(z)$ belongs to the convex hull of $M_z$
for all $z\in S$ is $k=\hat f$.
\endproclaim

\noindent{\bf Proof:}  
The assumptions are different from Theorem 4 in that the matrix
$(\alpha^i_j)$ is only defined for $i=1$ to $d-1$ and we require that
$d-1\leq \ell$ and $m\geq 2\ell$.  Since from (2) the $\partial_w q_i$
are pointwise linearly independent on $M$ ($i=1$ to $d-1$), the matrix
$(\alpha^i_j(z))$ has maximal rank for all $z\in S$, so for all
$z\in\overline D$.  By Theorem 2.2 of [SW], the matrix $(\alpha^i_j)$
can be extended to be defined for $i=1$ to $m$ as in Theorem 4, so
that the determinant of $(\alpha^i_j)$ is nonzero on $\overline D$.
Corollary $4$ then follows immediately from Theorem 4.
$\hfill\square$

We now present some examples.  Let $B_n$ be the open unit ball in
$\Bbb C^n$.

\proclaim{Example 1}  \endproclaim Let $g=(g_1,g_2,...,g_m)$ and
$k=(k_1,k_2,...,k_m)$ be $\Bbb C^m$-valued mappings analytic in a
neighborhood of the closed unit ball $\overline B_\ell$ in $\Bbb
C^\ell$.  Suppose that for each $i=1,2,...,m$, $g_i$ is never zero on
$\overline B_\ell$.  In (1), let us suppose that $S$ is the unit
sphere in $\Bbb C^\ell$, so $c=1$ and $p_1(z) =\|z\|_\ell^2-1$, where
$\|\cdot\|_n$ denotes Euclidean length in $\Bbb C^n$, $n\geq 1$.  In
(2), let $d=1$ and let $q_1(z,w)=\sum_{i=1}^m|g_i(z)w_i-k_i(z)|^2-1$,
so $M=\{(z,w)\in S\times\Bbb
C^m:\sum_{i=1}^m|g_i(z)w_i-k_i(z)|^2-1=0\}$.  We claim that the
mappings determined by Corollary $1$ are all mappings of the form
$$\eqalign{\overline B_
\ell&\rightarrow \Bbb C^m\cr z&\mapsto f(z)\equiv({k_1(z)+a_1\over g_1(z)},
 {k_2(z)+a_2\over g_2(z)},...,{k_m(z)+a_m\over g_m(z)})\cr}$$ where
$a=(a_1,a_2,...,a_m)\in\Bbb C^m$ is a constant of modulus $1$.  The
graphs of these maps in $S\times\Bbb C^m$ clearly foliate $M$.
Furthermore ${\partial q_1\over \partial
w_i}(z,w)=g_i(z)(\overline{g_i(z)w_i-k_i(z)})$ for $i=1,2,...,m$, so
${\partial q_1\over \partial w_i}(z,f(z))=\overline a_i g_i(z)$ for
all such $i$.  Since $\overline a_i$ is constant, it is CR on $S$ and
if we let $C(z)\equiv 1$, then (III) of Theorem 1 is satisfied.

Now suppose that $\ell=m$, $g_i=1$ and $h_i=0$ for $i=1$ to $m$.  Let
$T=(T_1,T_2,...,T_\ell)$ be an automorphism of $\overline B_\ell$.
Then the graph of $T$ over $\partial B_\ell$ in $\partial\overline
B_\ell\times\overline B_\ell$ lies in $M$.  However, such a graph is
not in the boundary of the polynomial hull of $M$ (which is $\overline
B_\ell\times
\overline B_\ell$) and hence by Corollary 3 these graphs do not
arise from elements of ${\Cal F}$.  Also, ${\partial q_1\over\partial
w_i}(z,T(z))=\overline {T_i(z)}$.  Given any $z^0\in\partial B_\ell$,
there is no nonzero $C(z)$ defined on $\partial B_\ell$ such that for
all $i$, $C(z)\overline{T_i(z)}$ is CR in $z$ near $z^0$; if there
were, then
$\sum_{i=1}^mC(z)\overline{T_i(z)}T_i(z)=\sum_{i=1}^mC(z)|T_i(z)|^2=C(z)$
would be CR near $z^0$.  Thus for all $i$, $\overline {T_i(z)}$ would
be CR in $z$ near $z^0$, which is impossible since this would imply
that the derivative of $T$ on the sphere is degenerate near $z^0$.

$\hfill\square$

\proclaim{Example 2} \endproclaim  Let $g=(g_1,g_2)$ be a $\Bbb C^2$-
valued analytic mapping defined in a neighborhood of $\overline
B_\ell$.  We suppose that $S=\partial B_\ell$ with
$p_1(z)=\|z\|_\ell^2-1$.  As for $M$, we let $d=2$ in (2), so that $M$
is defined in $S\times\Bbb C^2$ by two real valued functions $q_1,q_2$
defined as follows.  Let $h=(h_1,h_2)$ be a $\Bbb C^2$-valued mapping
analytic in a neighborhood of $\overline B_\ell$ which is never zero
on $\overline B_\ell$.  For $(z,w)\in S\times\Bbb C^2$ we let
$q_1(z,w)=\|w-g(z)\|_2^2-\|h(z)\|_2^2$ and let $q_2(z,w)=
2\Re(\sum_{i=1}^2h_i(z)(w_i-g_i(z)))$.  Consider the class of
functions
$$\eqalign{S&\rightarrow \Bbb C^m\cr z&\mapsto
g(z)+e^{i\theta}(-h_2(z),h_1(z))\cr}$$ where $\theta$ ranges over the
real numbers.  Then we claim that the graphs of these functions
foliate $M$ and in fact satisfy (III) of Theorem 1.  We let $\theta$
be an arbitrary real number and let $f(z)=g(z)+
e^{i\theta}(-h_2(z),h_1(z))$.  Then we calculate $q_1(z,f(z))=
\|g(z)+e^{i\theta}(-h_2(z),h_1(z))-g(z)\|_2^2-\|h(z)\|_2^2=
\|(-h_2(z),h_1(z))\|_2^2-\|h(z)\|_2^2=0$ for all $z\in S$ and
$q_2(z,f(z))=2\Re(\sum_{i=1}^2h_i(z)(f_i(z)-g_i(z)))=2\Re[(e^{i\theta})(-h_1(z)
h_2(z)+h_2(z)h_1(z))]=0$, as desired.  Furthermore, we calculate
$\partial_w q_1(z,w)=(\overline w_1-\overline g_1(z))dw_1+(\overline
w_2-\overline g_2(z))dw_2$ and $\partial_w
q_2(z,w)=h_1(z)dw_1+h_2(z)dw_2$.  Then $\partial_w
q_1(z,f(z))\wedge\partial_w q_2(z,f(z))=(h_2(z)(\overline
{f_1(z)}-\overline {g_1(z)})-h_1(z)(\overline{f_2(z)}-\overline
{g_2(z)}))dw_1\wedge dw_2
=-e^{-i\theta}(|h_1(z)|^2+|h_2(z)|^2)dw_1\wedge
dw_2=-e^{-i\theta}\|h(z)\|_2^2 dw_1\wedge dw_2$.  Letting $C(z)\equiv
\|h(z)\|_2^2$ in (III) of Theorem 1, we find that the properties of
(III) are satisfied for $f$ and $C$.  (The key fact is that $h$
is never zero on the closed ball.)  $\hfill\square$



\Refs
\widestnumber\key{BER}

\ref \key Al \by Alexander, H. \pages 201--212
\paper Polynomial hulls of graphs
\yr 1991 \vol 147, no. 2
\jour Pacific J. Math.
\endref

\ref \key AR \by Ahern, Patrick and Walter Rudin \pages 1--27
\paper Hulls of $3$-spheres in $\Bbb C\sp 3$
\yr 1992 \vol 137
\jour Contemp. Math. 
\endref

\ref \key B1 \by Begehr, Heinrich \pages 407--425
\paper Boundary value problems in $\Bbb C$ and $\Bbb C\sp n$
\yr 1997 \vol 22, no. 2
\jour Acta Math. Vietnam.
\endref

\ref \key B2 \by Begehr, Heinrich 
\book Complex analytic methods for partial differential equations
\publ World Scientific \publaddr River Edge, New Jersey \yr 1994
\endref

\ref \key B3 \by Begehr, Heinrich \pages 59--84
\paper Riemann-Hilbert boundary value problems in $\Bbb C\sp n$
\yr 1999 \vol 2
\jour Int. Soc. Anal. Appl. Comput.
\endref

\ref \key BD \by Begehr, Heinrich and Abduhamid Dzhuraev 
\book An introduction to several complex variables and partial 
differential equations
\publ Longman, Harlow \publaddr Harlow Essex \yr 1997
\endref

\ref \key BER \by Baouendi, M. Salah, Peter Ebenfelt and Linda Preiss
Rothschild
\book Real Submanifolds in Complex Space and Their Mappings
\publ Princeton University Press \publaddr Princeton, New Jersey \yr 1999
\endref

\ref \key Bo \by Boggess, Albert
\book CR Manifolds and the Tangential Cauchy-Riemann Complex
\publ CRC Press \publaddr Boca Raton, Florida \yr 1991
\endref

\ref \key \v Ce1 \by \v Cerne, Miran 
\paper Analytic discs in the polynomial hull of a disc fibration
over the sphere \toappear
\endref

\ref \key \v Ce2 \by \v Cerne, Miran 
\paper Maximal plurisubharmonic functions and the polynomial hull
of a completely circled fibration \toappear
\endref

\ref \key G \by Globevnik, Josip \pages 287--316
\paper Perturbation by analytic discs along maximal real submanifolds of
$\C^N$
\yr 1994 \vol 217, no. 2
\jour Math. Z.
\endref

\ref \key D1  \by Dzhuraev, Abduhamid \pages 271--295
\paper On linear boundary value problems in the unit ball of $\Bbb C\sp n$
\yr 1996   \vol 3, no. 2
\jour J. Math. Sci. Univ. Tokyo
\endref

\ref \key D2 \by Dzhuraev, Abduhamid \pages 287--303
\paper On Riemann-Hilbert boundary problem in several complex variables
\yr 1996  \vol 29, no. 4
\jour Complex Variables Theory Appl. 
\endref

\ref \key Fo \by Forstneri\v c, Franc \pages 869--889
\paper Polynomial Hulls of Sets Fibered Over the Circle
\yr 1988 \vol 37
\jour Indiana Univ. Math. J.
\endref

\ref \key Fr1 \by Freeman, Michael \pages 319--352
\paper Local biholomorphic straightening of real submanifolds
\yr 1977 \vol 106, no. 2
\jour Ann. Math. (2)
\endref

\ref \key Fr2 \by Freeman, Michael \pages 1--30
\paper Local complex foliation of real submanifolds
\yr 1974 \vol 209
\jour Math. Ann.
\endref

\ref \key Fr3 \by Freeman, Michael \pages 141--147
\paper Real submanifolds with degenerate Levi form
\inbook Several Complex Variables, Part 1
\publ American Mathematical Society \publaddr Providence, Rhode Island
\yr 1977
\endref

\ref \key Fr4 \by Freeman, Michael \pages 369--370
\paper The Levi form and local complex foliations
\yr 1976 \vol 57, no. 2
\jour Proc. Amer. Math. Soc.
\endref

\ref \key G \by Garnett, John B.
\book Bounded Analytic Functions
\publ Academic Press \publaddr New York \yr 1981
\endref

\ref \key Ha \by Han, Zu Hong \pages 701--708
\paper Foliations on CR manifold
\yr 1992 \vol 35, no. 6
\jour Sci. China Ser. A
\endref

\ref \key He \by Henkin, G.M. \pages 611--632
\paper Integral representation of functions which are holomorphic in 
strictly pseudoconvex regions, and some applications
\yr 1969 \vol 78
\jour Mat. Sb. (N.S.)
\endref

\ref \key H\" o \by H\" ormander, Lars 
\book Notions of Convexity
\publ Birkh\" auser \publaddr Boston
\yr 1994
\endref

\ref \key Hu \by Huang, Sha \pages 382--388
\paper A nonlinear boundary value problem for analytic functions of 
several complex variables 
\yr 1997 \vol 17, no. 4
\jour Acta Math. Sci.
\endref

\ref \key HMa \by Helton, J. William and Donald E. Marshall \pages 157--184
\paper Frequency domain design and analytic selections
\yr 1990 \vol 39, no. 1
\jour Indiana Univ. Math. J.
\endref

\ref \key HMe \by Helton, J. William and Orlando Merino \pages 285--287
\paper A fibered polynomial hull without an analytic selection
\yr 1994 \vol 41, no. 2
\jour Michigan Math. J.
\endref

\ref \key Ki \by Kiselman, Christer O. \pages 1--10
\paper A differential inequality characterizing weak lineal convexity
\yr 1998 \vol 311
\jour Math. Ann.
\endref

\ref \key Kr \by Kraut, Peter \pages 305--310
\paper Zu einem Satz von F. Sommer \" uber eine komplexanalytische
Bl\" atterung reeller Hyperfl\" achen im $\Bbb C\sp{n}$
\yr 1967 \vol 174 
\jour Math. Ann.
\endref

\ref \key Li \by Lieb, Ingo \pages 56--60
\paper Ein Approximationssatz auf streng pseudokonvexen Gebieten
\yr 1969 \vol 184
\jour Math. Ann.
\endref

\ref \key Ma \by Massey, William S. 
\book A Basic Course in Algebraic Topology
\publ Springer-Verlag \publaddr New York \yr 1991
\endref

\ref \key R \by Riemann, Bernhard \pages 3--47
\paper Grundlagen f\" ur eine allgemeine Theorie der Functionen
einer ver\" anderlichen complexen Gr\" osse
\inbook Gesammelte mathematische Werke und wissenschaftlicher Nachla\ss
\publ R. Dedekind and H. Weber \publaddr Leipzig
\yr 1876
\endref

\ref \key Ra \by Range, R. Michael 
\paper Extension phenomena in multidimensional complex analysis: 
correction of the historical record
\jour Math. Intelligencer \toappear
\endref

\ref \key S \by S\l odkowski, Zbigniew \pages 367--391
\paper Polynomial Hulls in {\bf C}$^2$ and Quasicircles
\yr 1989 \vol XVI 
\jour Ann. Scuola Norm. Sup. Pisa Cl. Sci. (4)
\endref

\ref \key Sh1 \by Shnirelman, A. I. \pages 257--258
\paper Degree of a quasiruled mapping, and the
nonlinear Hilbert problem
\yr 1972 \vol 27, no. 5
\jour Uspehi Mat. Nauk
\endref

\ref \key Sh2 \by Shnirelman, A. I. \pages 366--389, 533
\paper The degree of a quasiruled mapping, and the
nonlinear Hilbert problem
\yr 1972 \vol 89 (131)
\jour Mat. Sb. (N.S.)
\endref

\ref \key So1 \by Sommer, Friedrich \pages 111--133
\paper Komplex-analytische Bl\" atterung reeller Mannigfaltigkeiten
im $\Bbb C\sp{n}$
\yr 1958 \vol 136
\jour Math. Ann.
\endref

\ref \key So2 \by Sommer, Friedrich \pages 392--411
\paper Komplex-analytische Bl\" atterung reeler Hyperfl\" achen im
$\Bbb C\sp{n}$
\yr 1959 \vol 137
\jour Math. Ann.
\endref

\ref \key SW \by Sibony, N. and J. Wermer \pages 103--114
\paper Generators for $A(\Omega )$
\yr 1974 \vol 194
\jour Trans. Amer. Math. Soc.
\endref

\ref \key Tu \by Tumanov, A.E. \pages 385--398
\paper Extensions of CR-functions into a wedge
\yr 1991 \vol 70, no. 2
\jour Math. USSR-Sb.
\endref

\ref \key V \by Vekua, I.N. 
\book Generalized Analytic Functions
\publ Pergamon Press \publaddr London \yr 1962
\endref

\ref \key Wa \by Wang, Li Ping \pages 60--63
\paper Connective boundary value problems for analytic functions of
several complex variables in crisscross unbounded polycylindrical domains
\yr 1996 \vol 19, no. 2
\jour J. Sichuan Normal Univ.
\endref
 
\ref \key We1 \by Wegert, E. \pages 322--334
\paper Boundary value problems and best approximation by
holomorphic functions
\yr 1990 \vol 61, no. 3
\jour J. Approx Theory
\endref

\ref \key We2 \by Wegert, E. \pages 233--256
\paper Boundary value problems and extremal problems for holomorphic
functions
\yr 1989 \vol 11, no. 3-4
\jour Complex Variables Theory Appl.
\endref

\ref \key We3 \by Wegert, Elias
\book Nonlinear boundary value problems for holomorphic functions
and singular integral equations
\publ Akademie Verlag \publaddr Berlin \yr 1992
\endref

\ref \key We4 \by Wegert, E. \pages 583--615
\paper Nonlinear Riemann-Hilbert problems - history and perspectives
\yr 1999 \vol 11
\jour Ser. Approx. Decompos.
\endref

\ref \key We5 \by Wegert, E. \pages 307--313
\paper Nonlinear Riemann-Hilbert problems with unbounded 
restriction curves
\yr 1994 \vol 170
\jour Math. Nachr.
\endref

\ref \key Wa \by Warner, Frank W.
\book Foundations of Differentiable Manifolds and Lie Groups
\publ Springer-Verlag \publaddr New York
\endref

\ref \key Wh \by Whittlesey, Marshall A. \pages 677-701
\paper Polynomial hulls and $H^\infty$ control for a hypoconvex constraint.
\yr 2000 \vol 317
\jour Math. Ann.
\endref

\endRefs
\enddocument